\documentstyle[12pt,twoside]{amsart}

\title
{Arcs, valuations and the Nash map}
\author{Shihoko Ishii} 
\address{Department of Mathematics, Tokyo Institute of
Technology, Oh-Okayama, Meguro, Tokyo, Japan
\newline
e-mail : shihoko@@math.titech.ac.jp}

\newcommand{\bC}{{\Bbb C}}

\newcommand{\bZ}{{\Bbb Z}}

\newcommand{\bR}{{\Bbb R}}
\newcommand{\bN}{{\Bbb N}}
\newcommand{\bA}{{\Bbb A}}

\newcommand{\bm}{{\bf m}}
\newcommand{\codim}{\operatorname{codim}}
\newcommand{\Spec}{\operatorname{Spec}}

\newcommand{\Hom}{\operatorname{Hom}}
\newcommand{\Sing}{\operatorname{Sing}}
\newcommand{\ord}{\operatorname{ord}}
\newcommand{\val}{\operatorname{val}}
\newcommand{\trdeg}{\operatorname{trdeg}}
\newcommand{\st}{{\Spec \bC[[t]]}}

\newcommand{\sT}{{\Spec K[[t]]}}

\let \cedilla =\c
\renewcommand{\c}[0]{{\mathbb C}}  

\renewcommand{\o}[0]{{\mathcal O}} 

\newcommand{\spec}[0]{\operatorname{Spec}}
\newcommand{\sing}[0]{\operatorname{Sing}}

\newcommand{\GL}{\operatorname{GL}}
\newcommand{\SL}{\operatorname{SL}}

\newcommand{\vol}{\operatorname{vol}}

\def\to {\longrightarrow}

\newtheorem{thm}{Theorem}[section]

\newtheorem{lem}[thm]{Lemma}
\newtheorem{cor}[thm]{Corollary}
\newtheorem{prop}[thm]{Proposition}

\theoremstyle{definition}
\newtheorem{defn}[thm]{Definition}

\newtheorem{say}[thm]{}
\newtheorem{exmp}[thm]{Example}

\theoremstyle{remark}

\begin{document}
\maketitle
\footnote{partially supported by Grand-In-Aid of Ministry of Science 
and Education in Japan}

\begin{abstract}
  This paper gives a map from the set of families of arcs on a variety 
  to the set of valuations on the rational function field of the variety
  We characterize a family of  arcs which corresponds to a divisorial 
  valuation 
  by this map.
  We can see that both the Nash map and a certain McKay correspondence  are 
  the restrictions of this map.
  This paper also gives 
  the affirmative answer to the Nash problem for 
  a non-normal variety in a certain category.
  As a corollary, we get the affirmative answer for a non-normal 
   toric variety.

\end{abstract}

\section{Introduction}

\noindent
  In \cite{nash}, Nash introduces the Nash map which associates 
  a family of arcs through the singularities on a variety 
  (this family is called a Nash component in this paper) to 
  an essential divisor over the variety. 
  In other word, Nash map is a correspondence between the set of 
  certain families of arcs and
  the set of certain  divisorial valuations.
  
  On the other hand, L. Ein, R. Lazarsfeld and M. Musta\cedilla{t}\v{a} (\cite{ELM}) 
   introduce a map from the set of   
  irreducible cylinders for a non-singular variety to 
  the set of  
  divisorial
  valuations.

  In this paper, 
  we  introduce a map from the set of fat 
  arcs to the set of valuations.
  Here, a fat arc is an arc which does not factor through any 
  proper closed subvarieties.
  This map is a generalization of Nash map and the  map by Ein, 
  Lazarsfeld and Musta\cedilla{t}\v{a}.
  We can see that some fat arcs correspond to divisorial valuations 
  and the others to non-divisorial valuations.
  Here, we determine the fat arcs which correspond to divisorial 
  valuations. 
  By this characterization we obtain many examples corresponding to 
  divisorial valuations including Nash components and 
  cylinders in the arc space of a 
  non-singular variety. 
  As a cylinder and a Nash component are of  infinite dimension, 
  one may have an impression that an arc corresponding to a divisorial valuation 
  should be of infinite dimension.
  But  our characterization gives many finite dimensional families of arcs
  which correspond to  divisorial valuations.
  Another example  is the arc determined by a conjugacy class of a finite 
  group \( G \) which gives the quotient variety \( X=\bC^n/G \) 
  (\cite{d-l}).
  The restriction of our map onto a subset of these arcs coincides with the ``McKay 
  correspondence'' constructed in \cite{i-r}.

  This paper also study the Nash problem which asks if the Nash map is 
  bijective. 
  This problem was posed in Nash's preprint in 1968 
  (This preprint was published later as \cite{nash}).
  Inspired by this preprint, many people studied the arc spaces 
  of singularities and divisors over the singular varieties 
  (see,  Bouvier  \cite{B}, Gonzalez-Sprinberg \cite{GL}, 
  Hickel \cite{H},  Lejeune-Jalabert \cite{L80}, \cite{L}, \cite{LR},
  Nobile \cite{nobile}, Reguera-Lopez \cite{RL}) 
  Then, affirmative answer for the Nash problem is obtained for a 
  minimal 2-dimensional singularity by Reguera-Lopez \cite{RL}.
  For non-minimal 2-dimensional singularities, we do not know
  the answer of the Nash problem even for a rational double point
  (Recently the author was informed that a French 
  mathematician proved the affirmative answer 
  for a rational double point). 
  Last year, for a normal toric variety of arbitrary dimension  
  the Nash problem is affirmatively answered  
  but is negatively answered in general in \cite{i-k}.
 Though there is a counter example for the Nash problem, 
 it is still an interesting problem to clarify in which category 
 the Nash problem is affirmatively answered.
 For example, this problem is still open for 2 and 3 dimensional 
 singularities 
 as the counter examples in \cite{i-k} are normal singularity of 
 dimension greater than or equal to 4.
 For non-normal singularities, nothing is known about the Nash problem.
 In this paper, we give the affirmative answer to the Nash problem 
 for a non-normal singularity in a certain category.
 As a corollary, we obtain that for a non-normal toric variety
  the Nash problem is affirmative.

  This paper is organized as follows:
  In the second section, we give the basic notions of fat arcs.
  The map from the set of fat arcs to the set of valuations is given 
  here.
  
  In the third section, we give a characterization of a fat arc 
  which 
  corresponds to  a divisorial valuation. 
  Some examples including a cylinder on a non-singular variety  are 
  shown. The ``McKay correspondence'' in \cite{i-r} also appears as an 
  example.
  
  In the fourth section, we give the basic notions of the arc space of
  a toric variety, which are used in the following section.
  
  In the fifth section, we define a pretoric variety and prove that 
  the Nash problem is affirmative for a pretoric variety.
  As a non-normal toric variety is a pretoric variety, 
  the Nash problem is affirmative for this.

  Throughout this paper, the base field  is the complex number field 
  \( \bC \).
   A variety is an irreducible reduced scheme of finite type over \( 
   \bC \).
  A valuation is always a discrete valuation and 
  a valuation ring is a discrete valuation ring.
   
   The author would like to thank Professors Mircea 
   Musta\cedilla{t}\v{a} and K-i. Watanabe for 
   the stimulating discussions and valuable comments.
   She is grateful to Professor Hironobu Maeda for  providing 
   with information on valuation rings.
   The author is also grateful to  the members of Singularity Seminar 
   at Nihon University 
   for useful suggestions and encouragement.
   

\section{ The valuation associated to a fat arc}   

\begin{defn}
  Let \( X \) be a scheme of finite type over \( \bC\)
and $K\supset \bC$ a field extension.
  A morphism  \(\alpha: \Spec K[[t]]\to X \) is called an {\it arc} of \( X \).
  We denote the closed point of \( \Spec K[[t]] \) by \( 0 \) and
  the generic point by \( \eta \).
  The image \( \alpha(0) \) of the closed point is called the {\it 
  center} of the arc \( \alpha \).
  The transcendence degree of \( K\) over \( \bC\) is called 
  the {\it  dimension} of an arc \( \alpha \).
  Denote the space of arcs of \( X \) by \( X_{\infty} \).
\end{defn}

  If \( X \) is an affine variety, 
  the space \( X_{\infty} \) of arcs of \( X \) 
  is a closed subscheme of \( \Spec \bC[x_{1},x_{2},\ldots.] \),
  where \( \bC[x_{1},x_{2},\ldots.] \) is the polynomial ring over \( \bC
  \) with countably infinite number of variables.
  For a scheme \( X \) of finite type over \( \bC\),
  the arc space \( X_{\infty} \)
  is characterized by the 
  following property:
\begin{prop}
\label{ft}
  Let \( X \) be a scheme of finite type over \( \bC\).
  Then
  \[ \Hom _{\bC}(Y, X_{\infty})\simeq\Hom _{\bC}(Y\widehat\times_{\Spec \bC}\st, X) \]
   for an arbitrary \( \bC\)-scheme \( Y \),
   where \( Y\widehat\times_{\Spec \bC}\st \) means the formal completion 
   of \( Y\times_{\Spec \bC}\st \) along the subscheme 
   \( Y\times _{\Spec \bC} \{0\} \).
\end{prop}

\begin{say}
  By thinking of the case \( Y=\Spec K \) for an extension field \( K \) 
  of \( \bC\), 
  we see that  \( K \)-valued points of \( X_{\infty} \) correspond to
  arcs \( \alpha: \sT\to X \) bijectively. 
  Based on this, we denote the \( K \)-valued point corresponding to 
  an arc \( \alpha: \sT\to X \) by the same symbol \( \alpha \).
  The canonical projection \( X_{\infty} \to X \), \( \alpha\mapsto 
  \alpha(0) \) is  denoted by \( \pi_{X} \). 
  If there is no risk of confusion, we write   just \( \pi \).
 
  A morphism \( \varphi:Y \to X \) of varieties induces a canonical 
  morphism \( \varphi_{\infty}:Y_{\infty} \to X_{\infty} \), \( 
  \alpha \mapsto \varphi\circ\alpha \).

  The concept ``thin'' in the following is first introduced in 
  \cite{ELM}. 
\end{say}

\begin{defn}
   Let \( X \) be a variety over \( \bC\). 
   We say an arc \( \alpha:\sT  \to  X \) is {\it thin} 
   if \( \alpha  \) factors through a proper closed subvariety of \( X \).
  An arc which is not thin is called a {\it fat arc}.

  An irreducible subset \( C \) in \( X_{\infty} \)  is called a {\it 
  thin set} if \( C \) is contained in \( Z_{\infty} \) for a proper 
  closed subvariety \( Z\subset X \).
  An irreducible subset  in \( X_{\infty} \) which is not thin  is called a {\it 
  fat set}.
  
  In case an irreducible subset \( C \) has the generic point \( \gamma\in C 
  \) (i.e., the closure \( \overline{\gamma}  \) contains \( C \)),
  \( C \) is a fat set if and only if \( \gamma \) is a fat arc.

 \end{defn} 
\begin{prop}
\label{fatbasic}
  Let \( X \) be a variety over \( \bC\) and \( \alpha:\sT\to X \) an 
  arc. 
  Then, the following hold:
\begin{enumerate}
\item[(i)]
  \( \alpha \) is a fat arc if and only if
  the ring homomorphism \( \alpha^*: \o_{X, \alpha(0)}\to K[[t]] \) 
  induced from \( \alpha \) is injective;
\item[(ii)]
  Assume that \( \alpha \) is fat.
  For an arbitrary proper birational morphism  
  \( \varphi:Y\to X \),  \( \alpha \) is lifted to \( Y \).
\end{enumerate}  
\end{prop}

\begin{pf}
  As the problems are local, we may assume that \( X=\Spec A \) for a \( \bC
  \)-algebra \( A \).

(i).
 An arc \( \alpha:\sT\to X \) does not factor through any proper 
  closed subvariety of \( X \) if and only if  
  the generic point \( \eta \) of \( \sT \) is mapped to the generic 
  point of \( X \).
  This is equivalent to the injectivity of \(\alpha^*: A \to K[[t]] \)
  and also equivalent to the injectivity of \( \alpha^*: \o_{X, \alpha(0)}\to K[[t]]  \).

  (ii).
   Let \( \varphi:Y\to X \) be a proper birational morphism. 
  If the image \( \alpha(\eta) \) is the generic point of \( X \), 
  it is in the open subset on which \( \varphi \) is isomorphic. 
  Therefore, the restriction \( \Spec K((t))\to X \) of \( \alpha \) 
  is lifted to 
  \( Y \). 
  Then, by the criterion of properness, \( \alpha \) is lifted to \( Y \).      
\end{pf} 

\begin{defn}
  Let \( \alpha:\sT\to X \) be a fat arc of a variety \( X \) and 
  \( \alpha^*: \o_{X, \alpha(0)}\to K[[t]] \) the  ring homomorphism 
  induced from \( \alpha \).
  By \ref{fatbasic}, (i), \( \alpha^* \) is extended to the injective 
  homomorphism of fields \( \alpha^*:K(X)\to K((t)) \),
  where \( K(X) \) is the rational function field of \( X \).
  Define a function \(v_{\alpha}: K(X)\setminus\{0\}\to \bZ \) by
  \[ v_{\alpha}(f)=\ord_{t} \alpha^*(f). \] 
  Then, \( v_{\alpha} \) is a valuation of \( K(X) \).
  We call it the {\it  valuation corresponding to } \( \alpha \).
\end{defn}

\begin{prop}[Upper semicontinuity]
  For a regular function \( f \) on a variety \( X \), the map \( X_{\infty}\to 
  \bZ_{\geq 0} \), \( \alpha \mapsto v_{\alpha}(f) \) is upper 
  semicontinuous, 
  i.e., for \( r\in \bZ_{\geq 0} \) the subset 
  \( U_{r}= \{\alpha \in X_{\infty} \mid v_{\alpha}(f) \leq r \}\) is 
  open.
\end{prop} 

\begin{pf}
  We may assume that \( X \) is an affine variety \( \Spec A \). 
  Let \( X_{\infty} \) be \( \Spec A_{\infty} \) and 
  let \( \Lambda^*:A \to A_{\infty}[[t]]  \) be the ring homomorphism
  corresponding to the universal arc 
  \( X_{\infty}\widehat{\times}\spec \bC[[t]]\to X \) of \( X \).
  If we  write \( \Lambda^*(f)= a_{0}+a_{1}t+a_{2}t^2+\ldots\) with
  \( a_{i}\in A_{\infty} \), then \( X_{\infty}\setminus U_{r} \) is 
  the closed subset defined by \( a_{0}=a_{1}=\ldots=a_{r}=0 \).  
\end{pf}

\begin{defn}
  Let  \( X \) be a variety. We say that \( D \) is a {\it divisor 
  over} \( X \)
  if there is a proper birational morphism \( Y\to X \) from a normal 
  variety \( Y \) to \( X \) and \( D  \) is an irreducible divisor on \( Y \). 
\end{defn}

\begin{defn}
  A fat arc \( \alpha \) is called a {\it divisorial arc} if \( 
  v_{\alpha}=q \val_{D} \), 
  where \( q\in \bN \) and \( \val_{D} \) is the valuation defined by 
  an irreducible divisor \( D \) over \( X \).
  Assume that a fat set has the generic point.
 Then the fat set is called a {\it divisorial set} if the generic point is 
  a divisorial arc.
\end{defn}

\begin{prop}
\label{divisorial}
  For a fat arc \( \alpha \) of a variety \( X \) of dimension \( n \), 
  the following are equivalent:
\begin{enumerate}
\item[(i)]
  \( \alpha \) is divisorial;
\item[(ii)]
   There is a proper birational morphism \( Y\to X \) from a normal 
   variety \( Y \) such that the center \( \tilde\alpha(0) \) 
   of the lifting \( \tilde \alpha \) 
   of \( \alpha \) to \( Y \) is the generic point of a divisor on \( Y \);
\item[(iii)]
  Identifying \( K(X) \) with the subfield of \( K((t)) \) by \( \alpha^* 
  \),
  we denote \( K[[t]]\cap K(X) \) by \( R \) and  \( (t)\cap K(X) \) by \( 
  \frak{m} \);  
  Then the transcendence degree \( \trdeg_{\bC} R/{\frak{m}}=n-1 \).
\end{enumerate}  
\end{prop}

\begin{pf}
  Assume (i), then \( v_{\alpha}=q\val_{D} \) for \( q\in \bN \) and 
  a divisor \( D  \) 
  on \( Y \), where \( Y \) is normal and there is a proper 
  birational morphism \( Y\to X \).
  Let \( \tilde \alpha \) is the lifting of \( \alpha \) to \( Y \).
  Let \( \delta \in Y \) be the generic point of \( D \), 
  then, by \( v_{\alpha}=q\val_{D} \), we have that 
  \[ \ord_{t}\alpha^*f\geq 0 \ \operatorname{for\ every\ } f\in 
  \o_{Y,\delta}\operatorname{\ and}\]
  \[  \ord_{t}\alpha^*f > 0 \ \operatorname{for\ every\ } f\in 
  {\frak{m}}_{Y,\delta}, \]
  which means that the field homomorphism 
  \( \alpha^*=\tilde \alpha^*:K(X)=K(Y)\to K((t)) \) gives a local 
  homomorphism \( \o_{Y,\delta}\to K[[t]] \). 
  Hence, \( \tilde\alpha(0)=\delta \) as required in (ii).
  
  Assume (ii), then \( \tilde\alpha \) induces an injective local 
  homomorphism \( \o_{Y,\delta}\to K[[t]] \), where \( \delta \) is 
  the generic point of the divisor over \( X \).
  Then \( R/{\frak{m}}\supset \o_{Y,\delta}/{\frak{m}}_{Y,\delta} \), 
  where the transcendence degree of the right hand side over \( \bC\) is \( n-1 \).
  Therefore, \( \trdeg_{\bC}R/{\frak{m}} \) must be \( n-1 \). 
  
  Assume (iii), then by Zariski's local uniformizing theorem
  (\cite{zariski}, see also \cite[II, \S 14, Theorem 31]{z-s}) 
  there is a divisor \( D \) over \( X \) such that 
  \( R=\o_{Y,\delta} \), where \( \delta \) is the generic point of \( 
  D \). 
  Since the valuation rings of \( v_{\alpha} \) and \( \val_{D}  \) are the same,
  it follows that \( v_{\alpha}=q\val_{D} \) for some \( q\in \bN \), 
  which implies (i).
\end{pf}

\begin{prop}
\label{existence}
   For every divisor \( D \) over a variety \( X \) and every \( q\in 
   \bN \), 
   there is a fat arc \( \alpha\in X_{\infty}  \) of \( X \) such that \( 
   v_{\alpha}=q\val_{D} \).
\end{prop}

\begin{pf}
  Let \( \varphi:Y\to X \) be a proper birational morphism such that \( D \) is 
  a divisor on a normal variety \( Y \).
  Let \( \delta \) be the generic point of \( D \).
  Replacing \( X \) by an affine open neighborhood of \( \varphi (\delta) 
  \), 
  we may assume that \( X \) is affine variety \( \Spec A \).
  Then, there are injections:
  \[ A \hookrightarrow \o_{Y,\delta}\hookrightarrow 
  {\widehat{\o_{Y,\delta}}} \simeq K[[t]],\]
 
  \noindent
  where \( \  \widehat{ }\     \) is the completion by the maximal ideal and 
  \( K \) is the residue field of the  local ring \( \o_{Y,\delta} \). 
  Composing these maps and the homomorphism \( K[[t]]\to K[[t']] \), 
  \( t \mapsto {t'}^q \), 
  we obtain an arc \( \alpha':\Spec K[[t']]\to X=\Spec A \).
  It is easy to see that \( v_{\alpha'}=q\val_{D} \). 
  Then, just take an image of \( \alpha':\Spec K \to X_{\infty} \)
  as  a required \( \alpha \).

\end{pf}

\begin{exmp} 
  Let \( X=\bA^2_{\bC} \) and \( \alpha:\st \to X \) the arc defined by 
  the ring homomorphism \( \alpha^*: \bC[x,y] \to \bC[[t]] \), 
  \( x \mapsto t \), \( y \mapsto e^t-1=t+t^2/2!+t^3/3!+\ldots \).
  Then \( \alpha \) is an arc of \( X \) with the center at 0.
  As \( \alpha^* \) is injective, the arc \( \alpha \) is fat by 
Proposition  \ref{fatbasic}.
  But it is not divisorial, 
  because \( R/{{\frak{m}}} \) defined in (iii) of Proposition \ref{divisorial} 
  is contained in \( \bC[[t]]/(t)=\bC \).
\end{exmp}

\begin{exmp}[A cylinder on a non-singular variety \cite{ELM}]
\label{cylinder1}
   Let \( X \) be a non-singular variety and \( C\subset X_{\infty} 
   \) an irreducible cylinder. 
  The paper \cite{ELM} defines a valuation \( \val_{C} \) of \( K(X) \) 
  corresponding 
  to \( C \) and proves that 
  \( \val_{C} \) is a divisorial valuation.
  It is easy to see that this valuation \( \val_{C} \) is the same as our
  valuation \( v_{\gamma} \), 
  where \( \gamma \) is the generic point of \( C \).
  In the next section, we will see another proof for the fact that 
  \( v_{\gamma} \) is divisorial.     
\end{exmp}

\begin{exmp}
\label{maximal}
  Let \( \varphi:Y \to X \) be a resolution of the singularities of \( X 
  \).
    Let \( E \) be an irreducible divisor on \( Y \).
  Let \( \beta  \) be the generic point of an irreducible 
  closed subscheme \( \pi_{Y}^{-1}(E)\subset Y_{\infty} \).
  Then, \( \beta \) is the lifting of \( 
  \alpha=\varphi_{\infty}(\beta) \) to \( Y \) and \( \beta(0) \) is the generic 
  point of \( E \).
  Therefore, by Proposition \ref{divisorial}, \( \alpha \in X_{\infty}\) is 
  a divisorial arc. 
  Actually it follows \[  v_{\alpha}=\val_{E} .
 \]  To see this, let \( v_{\alpha}=q\val_{E} \) and \( \gamma\in 
 X_{\infty}
  \) be an arc such that \( v_{\gamma}= \val_{E} \) (Proposition 
  \ref{existence}).
  Then, the lifting \( \tilde \gamma \) of \( \gamma  \) to \( Y \) 
  should have center \( \tilde\gamma(0) \) at the generic point of \( E \).
  Therefore \( \tilde \gamma \in \pi_{Y}^{-1}(E) \).
  As \( \tilde \gamma \) is contained in the closure of \( \beta \),
  it follows that 
  \(v_{\alpha}(f)= v_{\beta}(f)\leq v_{\tilde \gamma}(f)=v_{\gamma}(f) \) 
  for every \( f\in 
  \o_{Y,e} \) by upper semicontinuity, 
  where \( e \) is the generic point of \( E \).
  This yields \( q=1 \).
\end{exmp}

\begin{exmp}[Nash components, a special case of Example \ref{maximal}]
\label{nash}
  Let \( X \) be a variety and \( \Sing X \) the singular locus of \( 
  X \). 
  An irreducible component  \( C \) of \( \pi_{X}^{-1}(\Sing X) \) is called a 
  Nash component if \( C  \) is not contained in \( (\Sing X)_{\infty} 
  \). 
  (In \cite{i-k} a Nash component is called a ``good component''.)
  By \cite[Lemma 2.12]{i-k} 
   every irreducible component of \( \pi_{X}^{-1}(\Sing 
  X) \) 
  is a Nash component if the base field is of characteristic zero.
  Noting that our base field is \( \bC \), let \( \pi_{X}^{-1}(\Sing 
  X) = \bigcup_{i}C_{i} \) be the decomposition into Nash components.
  Let \( \varphi:Y \to X \) be a resolution such that \( \varphi \) 
  is isomorphic outside of \( \Sing X \) and \( \varphi^{-1}(\Sing X) \) 
  is a divisor. 
  Let \( \varphi^{-1}(\Sing X)=\bigcup_{j} E_{j} \) be the decomposition 
  into irreducible components.
  Then \( \varphi_{\infty}: \bigcup_{j} \pi_{Y}^{-1}(E_{j}) \to 
  \bigcup_{i}C_{i} \) is bijective outside of \( (\Sing X)_{\infty} \).
  Hence, for each \( C_{i} \) there is unique \( E_{j} \) such that 
  \( \pi_{Y}^{-1}(E_{j}) \) is dominant to \( C_{i} \).
  Therefore, the generic point \( \beta \) of \( \pi_{Y}^{-1}(E_{j}) \) 
  is mapped to the generic point \( \alpha \) of \( C_{i} \) by the 
  morphism \( \varphi_{\infty} \). 
  In \cite{nash} Nash proved that this  \( E_{j} \) is an essential divisor 
  over \( X \) 
  (for the proof see also \cite[Theorem 2.15]{i-k}).
  This map
  \begin{center}
   \( {\cal N}: \{ \) Nash components \( \}\to \{ \) essential 
  divisors over \( X \}\), \( C_{i } \mapsto E_{j} \) 
  \end{center}
  is called Nash map and Nash problem is the problem if this map is bijective. 
  By the discussion in Example 
  \ref{maximal}, it follows that \( v_{\alpha}=\val_{E_{j}} \).
  Hence,  a Nash component is divisorial and our correspondence 
  between fat arcs and the valuations gives the Nash map.        
\end{exmp}

\begin{exmp}
\label{toric}
  Here, we use the notation and terminologies of \cite{fulton}.
  Let  $M$  be the free abelian group  ${\bZ}^n$ $(n\geq 1)$
  and  $N$   its dual $\Hom_{\bZ}(M, {\bZ})$.
  We denote  $M\otimes _{\bZ}{\bR}$  and $N\otimes_{\bZ}{\bR}$  by
  ${M_{\bR}}$  and  $N_{\bR}$, respectively.
  The canonical pairing \( \langle\ , \ \rangle:
  N\times M \to \bZ \)  extends to  
  \( \langle\ \ , \ \rangle:
  N_{\bR}\times M_{\bR} \to \bR \).
  A cone in \( N_{\bR} \) generated by \( v_{1},\ldots,v_{n}\in N \) is denoted by 
  \( \langle v_{1},\ldots,v_{n}\rangle \).
  The group ring \( \bC[M] \) is generated by monomials \( x^{\bm} \) 
  (\( \bm \in M \)) over \( \bC \).
  Let \( X \) be an affine toric variety defined by a cone \( \sigma \) 
  in \( N \).
  In \cite{i}, for \( v\in \sigma\cap N \) we define 
 \[  T_{\infty}(v)= \{\alpha\in X_{\infty} \mid \alpha(\eta)\in T,\  \ord_{t} 
  \alpha^*(x^u)=\langle v, u \rangle \ {\operatorname{for}}\ u\in M \}, \]
   where \( T \) denotes the open orbit and also the torus acting on \( X \).
   Then,    \( T_{\infty}(v) \) is a divisorial fat set corresponding 
  to \( D_{v} \), where, for the maximal \( q\in \bN \) such that  \( 
  v/q\in N 
  \) and \( D_{v} \) means the divisor \( q 
  (\overline{orb \bR_{\geq 0}v )} \).
  This is proved as follows:
  The paper \cite{i} proves that \( T_{\infty}(v) \) is irreducible and 
  locally closed.
  Let \( \alpha \) be the generic point of \( T_{\infty}(v) \).
  Let \( \gamma\in X_{\infty} \) be a divisorial arc corresponding to 
  \( D_{v} \) (Proposition \ref{existence}).
  Then, by the definition of \( T_{\infty}(v) \), for \( u\in \sigma^{\vee}\cap M \)
  \[ v_{\alpha}(x^u)= \ord_{t}\alpha^*(x^u)=\langle v, u\rangle
  = \val_{D_{v}}(x^{u})=v_{\gamma}(x^u). \] 
For a general \( f\in \bC[ \sigma^{\vee}\cap M]\) we have that
\[   v_{\alpha}(f)= \ord_{t}\alpha^*(f)\geq \min_{x^u\in f}\langle v, u\rangle
  = \val_{D_{v}}(f)=v_{\gamma}(f). \]
    Here, \( {x^u\in f} \) means that the coefficient of \( x^u \) of \( 
  f \) is not zero.

  On the other hand, since \( \gamma  \) belongs to \(  T_{\infty}(v) 
  \) whose generic point is 
  \( \alpha \), it follows that 
  \( v_{\alpha}(f)\leq v_{\gamma}(f) \) 
  for every \( f\in \bC[\sigma^{\vee}\cap M] \) by upper semicontinuity.
  Therefore, \( v_{\alpha}=v_{\gamma}=\val_{D_{v}} \).  
\end{exmp}

  In particular for a primitive \( v\in \sigma \cap N \) we obtain 
  the following which will be used later in this paper.

\begin{lem}
\label{equal}
  Let \( X \) be an affine toric variety defined by a cone \( 
  \sigma  \) in \( N \) and \( \varphi:Y \to X \) be an 
  equivariant resolution of the singularities of \( X \).
  Let \( D \) be an irreducible invariant divisor on \( Y \) 
  corresponding to \( v\in \sigma \cap N \).
  Then the following closures coincide:
  \[  
  \overline{\varphi_{\infty}(\pi_{Y}^{-1}(D))}=\overline{T_{\infty}(v)}. \]
\end{lem}

\begin{pf}
  It is sufficient to prove that the generic point of each hand side 
  is contained in 
  the other hand side.
  Let \( \alpha \) be the generic point of \( 
  \varphi_{\infty}(\pi_{Y}^{-1}(D)) \).
  Then, by Example \ref{maximal}
  \[ \ord_{t}\alpha^*(f)=\val_{D}(f)=\min_{x^u\in f}\langle v,u 
  \rangle, \]
  for a regular function \( f \) on \( X \).
  Hence, in particular, \( \ord_{t}\alpha^*(x^u)=\langle v, u\rangle 
  \) for \( u\in \sigma^{\vee}\cap M \), which implies \( \alpha\in 
  T_{\infty}(v) \), 
  as \( \sigma^{\vee}\cap M \) generates the group \( M \).
  Conversely, if \( \alpha \) is the generic point of \( T_{\infty}(v) 
  \),
  \( \alpha \) is a divisorial arc corresponding to \( D \) by 
  \cite[Proposition 5.7]{i}.
  Let \( \tilde \alpha \) be the lifting of \( \alpha \) to \( Y \).
  Then, by Proposition \ref{divisorial}, \( \tilde\alpha(0) \) is the 
  generic point of \( D \).
  Therefore, \( \tilde\alpha\in \pi_{Y}^{-1}(D) \), which yields that
  \( \alpha=\varphi_{\infty}(\tilde\alpha)\in \varphi_{\infty}(\pi_{Y}^{-1}(D)) \).          
\end{pf}


\section{A characterization of a divisorial arc}

\noindent
  In this section, we characterize a divisorial arc.
  For this we start with a simple lemma.
  We note that this lemma follows immediately from Corollary 1 of 
  \cite[VI, \S 6]{z-s},
  when both \( R, R' \) in the statement of the lemma are valuation rings.
  
\begin{lem}
\label{algebraic}
  Let \( K\supset K' \) be an algebraic extension of fields. 
  Let \( R \) and \( R' \) be valuation rings in \( K \) and \( K' \), 
  respectively. Denote the maximal ideals of \( R \) and \( R' \)
  by \( {\frak{m}} \) and \( {\frak{m}}' \), respectively.
  Assume that \( R \) dominates \( R' \). 
  Then the field extension \( R/{\frak{m}}\supset  R'/{\frak{m}}'\)
  is algebraic. 
\end{lem}  

\begin{pf}
  Let \( v' \) be a valuation whose valuation ring is \( R' \).
  Let \( \overline{f} \) be an arbitrary element of \(  R/{\frak{m}} \) 
  and \( f \in R \)  an element corresponding to \( \overline{f} \).
  Then, there is an equality
  \[a_{n}f^n+a_{n-1}f^{n-1}+..+a_{0}=0,  \]
  with \( a_{j}\in R' \) (\( j=0,\ldots,n \)).
  Let \( v'(a_{i}) = \min_{j=0,\ldots,n}v'(a_{j}) \).
  Put \( b_{j}=a_{i}^{-1}a_{j} \) for \( j=0,..,n \).
  Note that \( b_{j}\in R' \) for all \( j \) and \( b_{i}=1 \).     
  Then we obtain the equality in \( R/{\frak{m}} \):
  \[ 
  \overline{b_{n}}\overline{f^n}+\overline{b_{n-1}}\overline{f^{n-1}}+
  \ldots+\overline{b_{0}}=0, \]
  which is an algebraic equation of \( \overline{f} \)  over  \( R'/{\frak{m}}' \). 
\end{pf}

\begin{lem}
\label{commutative}
  Let \( X \) and \( X' \) be varieties of the same dimension with a 
  dominant morphism \( \varphi:X\to X' \). 
  Let \( K\supset K' \) be a field extension.
  Assume there exists a commutative diagram of arcs
  \[ \begin{array}{ccc}
   \Spec K[[t]]& \stackrel{\alpha}\longrightarrow & X\\
   \downarrow & & \  \downarrow \varphi\\
   \Spec K'[[t]]& \stackrel{\alpha'}\longrightarrow & X'.\\
   \end{array} \]
  Then, \( {\alpha} \) is divisorial if and only if \( {\alpha'} 
  \) is divisorial.
\end{lem}

\begin{pf}
  Let \( R=K[[t]]\cap K(X) \) and \( R'=K'[[t]]\cap K(X') \) and 
  let \( {\frak{m}} \) and \( {\frak{m'}} \) be the maximal ideals of 
  \( R \) and \( R' \), respectively.
  Then \( R \) dominates \( R' \).
  Since \( K(X) \) is algebraic over \( K(X') \), 
  the field extension \( R/{\frak{m}}\supset R'/{\frak{m'}} \) is 
  algebraic by Lemma \ref{algebraic}.
  Therefore, \( \trdeg_{\bC}R/{\frak{m}}=n-1 \) if and only if 
  \(  \trdeg_{\bC}R'/{\frak{m'}}=n-1 \). 
  Then, by Proposition \ref{divisorial} we obtain the assertion.
\end{pf}

\begin{lem}
\label{minimal}
  Let \( K \supset K'\) be a field extension and \( 
  f_{i}=\sum_{j}a_{ij}t^j \) (\( i=1,\ldots,s \)) elements of 
  \( K[[t]] \).
  Assume that for each \( i \) \(( i=1,\ldots,s )\) there is \( j_{i} \)
  such that \( a_{1 j_{1}},a_{2 j_{2}},\ldots, a_{s j_{s}} \) are 
  algebraically independent over \( K' \).
  Then, by changing the numbering of \( \{i\} \), 
  for each \( i \) we have \( l_{i} \)
  such that \( a_{1 l_{1}},a_{2 l_{2}},\ldots, a_{s l_{s}} \) are 
  algebraically independent over \( K' \) and the following hold:
\begin{enumerate}
\item[(1)]
\( a_{1j}\in \overline{K'} \) for every \( j< l_{1} \);
\item[(2)]
\( a_{2j}\in \overline{K'(a_{1l_{1}})} \) for every \( j< l_{2} \);
\newline
\ldots..
\item[(s)]
\( a_{sj}\in \overline{K'(a_{1 l_{1}},\ldots, a_{s-1 l_{s-1}})} \) 
for every \( j<l_{s} \),
where \( \overline{K'(*) } \) means the algebraic closure of \( K'(*) \) 
in \( K \).
\end{enumerate}  
\end{lem}

\begin{pf}
  First we define a partial order \( \leq \) in \( \bZ_{\geq 0}^s \) 
  as follows:
  \[ (b_{1},\ldots,b_{s})\leq  (b'_{1},\ldots,b'_{s}) \ \ 
  \operatorname{if}\  b_{i}\leq b'_{i} \ \ \operatorname{for\ 
  every\ } i .\]
  Let \( (l_{1},\ldots,l_{s}) \) be a minimal element in 
  \[M_s= \{(j_{1},\ldots,j_{s}) \mid  a_{1 j_{1}},a_{2 j_{2}},\ldots, a_{s j_{s}} 
 \]
 \[\operatorname{\  are\  
  algebraically\  independent\  over\ }   K' \}.    
 \]
  Let  \( A_{i}=\{a_{ij}\mid j<l_{i}\}  \).
  Then,  every \( a_{ij}\in A_{i} \) is an element of  
  \( \overline{K'(a_{1 l_{1}},\ldots, a_{sl_{s}})} \). 
  Indeed, if \( a_{ij}\in A_{i} \) was transcendental over
  \( K'(a_{1 l_{1}},\ldots, a_{sl_{s}})  
 \),
   then
  \( (l_{1},..,l_{i-1},j, l_{i+1},..,l_{s}) \) would be a smaller 
  element than 
  \( (l_{1},..,l_{i-1},l_{i}, l_{i+1},..,l_{s}) \) in \( M_{s} \), which is a contradiction to the minimality of
   \( (l_{1},\ldots,l_{s}) \).
   Let \( B_{i}\subset \{a_{1l_{1}},\ldots,a_{sl_{s}}\} \) be a 
   minimal subset such that 
   \( A_{i}\subset \overline{K'(B_{i})} \).
   We prove by induction on \( s \) that these \( l_{i} \)'s are required 
   in our lemma.
   
   (I) The case \( s=1 \). By the minimality of \( l_{1} \),
   every element \( a_{1j}\in A_{1} \) is algebraic over \( K' \).
   
   (II) The case \( s\geq 2 \). 
   Assume the assertion  for \( s-1 \).
  By changing the numbering of \( \{i \} \), 
  let \( \# B_{1} \) attain  \( \min_{i} \# B_{i} \).   
  We can see that \( (l_{2},\ldots,l_{s}) \) is a minimal element in 
  \[M_{s-1}= \{(j_{2},\ldots,j_{s}) \mid  a_{2 j_{2}},a_{3 j_{3}},\ldots, a_{s j_{s}} 
 \]
 \[\operatorname{\  are\  
  algebraically\  independent\  over\ }   K'(a_{1l_{1}}) \}.    
 \]
  Indeed, if \( (j_{2},\ldots,j_{s}) \) was a smaller element than 
  \( (l_{2},\ldots,l_{s}) \) in \( M_{s-1} \),
  then \( (l_{1},j_{2},\ldots,j_{s}) \) would be a smaller element than
  \( (l_{1},l_{2},\ldots,l_{s}) \) in \( M_{s} \),
  which is a contradiction to the definition of 
   \( (l_{1},l_{2},\ldots,l_{s}) \).

  Then by the hypothesis of induction, we obtain
\begin{enumerate}
\item[(2) ]
\( A_{2}\subset \overline{K'(a_{1l_{1}})} \),
\item[(3) ]
\( A_{3}\subset \overline{K'(a_{1l_{1}},a_{2l_{2}})} \),
\newline
\ldots..
\item[(s) ]
\( A_{s}\subset  \overline{K'(a_{1l_{1}},a_{2l_{2}},.., a_{s-1l_{s-1}})} \).  
\end{enumerate}
  By this, we can see that \( \#B_{2}\leq 1 \).
  Hence, by the minimality of \( \#B_{1} \), we have \( \# B_{1}\leq 1 \).
  If \( \# B_{1}=0 \), then the proof is over.

  Assume \( \# B_{1}=1 \) and induce a contradiction.
  Let \( I=\{i \mid \# B_{i}=1\} \) and assume \( I=\{1,2,..,r\} \)
  \( (r\leq s) \).
  
Letting \( B_{i} \) \( (i\in I) \) play the role of \( B_{1} \), 
we carry out the above discussion.
  Then, we have the inclusion \( A_{\sigma(i)}\subset 
  \overline{K'(a_{il_{i}})} \) corresponding to (2) for some 
  \( \sigma(i)\in I \).
  This map \( \sigma:I\to I \) is bijective.
  Indeed, if \( \sigma(i)=\sigma(i') \) for \( i\neq i' \), 
  it would follow
  \( A_{\sigma(i)}\subset \overline{K'(a_{il_{i}})}
  \cap \overline{K'(a_{i'l_{i'}})}=\overline{K'} \), 
  which contradicts to the minimality of \( \# B_{1}=1 \).
  For every \( i\in I \), take an element \( a_{ik_{i}}\in A_{i} \) 
  such that \( a_{ik_{i}}\in \overline{K'(a_{i'l_{i'}})}\setminus 
  \overline{K'} \),
  where \( i'=\sigma^{-1}(i) \).
  Then, \( a_{1k_{1}},.., a_{rk_{r}}, a_{r+1l_{r+1}},..,a_{sl_{s}} \) 
  are algebraically independent over \( K'  \) and 
  \( (k_{1},..,k_{r},l_{r+1},..,l_{s})<(l_{1},..,l_{r},l_{r+1},..,l_{s}) 
  \),
  which is a contradiction.    
\end{pf}

\begin{lem}
\label{extension}
  Let \( K   \) be a field.
  For every element \( f\in K[[t]] \) of order \( d > 0 \), 
  there exists an algebraic field extension \( L\supset K \) 
  and  an element \( t'\in L[[t]] \) such that  \( L[[t]]= L[[t']] \) 
  and \( f={t'}^d \) in \( L[[t']] \).
\end{lem}

\begin{pf}
  Let \( f= a_{d}t^d+ a_{d+1}t^{d+1}+\ldots \) (\( 
  a_{d},a_{d+1},\ldots\in K \)).
  Let \( L \) be the extension field of \( K \) by adding the roots of 
  an equation \( X^d-a_{d}=0 \). 
  Put \( t'=b_{1}t+b_{2}t^2+\ldots \).
  Then the equation \( a_{d}t^d+ a_{d+1}t^{d+1}+\ldots = 
  (b_{1}t+b_{2}t^2+\ldots)^d \) has a solution 
  \( b_{1}, b_{2},\ldots \) in \( L \).
\end{pf}

\begin{defn}
  The extension \( L[[t']]\supset K[[t]] \) in Lemma 
   \ref{extension} is called 
  a {\it canonical extension with respect to} \( f \).
\end{defn}

\begin{thm}
\label{thm}
  Let \( \alpha:\sT\to X \) be a fat arc of a variety \( X \) of 
  dimension \( n \) and
  \( e  = \alpha(0) \).
  Let \( \alpha^*: {\o}_{X,e}\to K[[t]] \) be the ring homomorphism 
  corresponding to \( \alpha \).
  Let the codimension of the closure \( \overline{\{e\}} \) be \( 
  r\geq 2 \).
  Then the following are equivalent:
\begin{enumerate}
\item[(i)]
  \( \alpha \) is divisorial;

\item[(ii)]
 There exist \( x_{1},\ldots,x_{r}\in {\frak{m}}_{X,e} \) such 
     that, for a canonical extension \( L[[t']]\supset K[[t]] \) with
     respect to \( \alpha^*(x_{r}) \), 
     we have \( \alpha^*(x_{i})= \sum_{j\geq k_{i}}a_{ij}{t'}^j\) 
     with  a set of coefficients \( 
     a_{1l_{1}},a_{2l_{2}},\ldots,a_{r-1l_{r-1}} \) 
     (\( l_{i}\geq k_{i} \)) which are algebraically 
     independent over \( K'=\o_{X,e}/{\frak{m}}_{X,e} \),
      where \( {\frak{m}}_{X,e} \) is the maximal ideal of \( \o_{X,e} 
     \) and we regard \( K' \) as a subfield in \( K \) by \( \alpha^* \).     
\end{enumerate}  
\end{thm}

\begin{pf}
   First we may assume that \( X \) is affine. 
  
   (i)\( \Rightarrow \)(ii). 
   By Noether's normalization lemma, 
   there is a finite dominant morphism from \( X \) to a non-singular 
   variety of dimension \( n \).
   Then, by Lemma \ref{commutative} we may assume that \( X=\bA^n_{\bC}
   =\Spec \bC[x_{1},..,x_{r},..,x_{n}] 
   \) and \( {\frak{m}}_{X,e} \) is generated by \( r \)-elements 
   \( x_{1},\ldots,x_{r} \).
   Then,  \( K'=\o_{X,e}/{\frak{m}}_{X,e}=\bC(x_{r+1},.,x_{n}) \).  
    Let \( L[[t']]\supset K[[t]] \) be a canonical extension with
     respect to \( \alpha^*(x_{r})={t'}^d \).
     Then, the arc \( \alpha':\Spec L[[{t'}]]\to X \) induced from 
     \( \alpha \) is lifted as 
    \[   \alpha': \Spec L[[t']]\to X'=\Spec \bC[x_{1},.,x_{r-1}, 
     x_{r}^{\frac{1}{d}},x_{r+1},.,x_{n}]  \] and 
     this is divisorial by Lemma \ref{commutative}.
  Hence, there are a divisor \( D \) over \( X' \) and a 
  local homomorphism \( {\alpha'}^*:\o_{Y,\delta}\to L[[t']] \),
  where \( Y\supset D \) and \( \delta \) is the generic point of \( D \).
  Since \( {\alpha'}^*(x_{r}^{\frac{1}{d}})=t' \), it follows that \( 
  {\frak{m}}_{Y,\delta}=(x_{r}^{\frac{1}{d}}) \).
  Write \( x_{i}={x_{r}^{\frac{k_{i}}{d}}}u_{i} \) for a unit 
  \( u_{i} \) in \( \o_{Y,\delta} \).
  Here we note that \( a_{ij} \)'s are coefficients of 
  \( {\alpha'}^*(u_{i}) \).
  Assume  \( a_{ij}\in \overline{K'} \) for all \( i,j \),
  where   \( \overline{K'}  \) is the algebraic closure of \( K' \) in 
  \( L \).
Then,  \(   L[[t']]\cap K'(u_{1},..,u_{r-1},  x_{r}^{\frac{1}{d}}) 
  \subset \overline{K'}[[t']] \) which yields 
 \[ \left(L[[t']]\cap K(Y)\right) /\left((t')\cap K(Y)\right)\subset 
 \overline{ K'}[[t']]/(t')=
  \overline{K'} . \]
  Therefore  \[ \trdeg_{\bC}(L[[t']]\cap K(Y) )/((t')\cap K(Y))= n-r < n-1, \] 
  which is a contradiction to Proposition \ref{divisorial}.
 Hence, there exists \( a_{ij} \) transcendental over \( K' \).
  Let it be \( a_{1l_{1}} \) and assume that \( a_{1j}\in \overline{K'} \)
  for \( j< l_{1} \).
  Let \( R_{1}:=L[[t']]\cap \overline{K'}(u_{1}, x_{r}^{\frac{1}{d}})\)
  and \( {\frak{m}}_{1} \) the maximal ideal.
  Then, \( L[[t']]/(t')\supset R_{1}/{\frak{m}}_{1}=\overline{K'}(a_{1l_{1}}) \).
  This is proved as follows:
  Let \[   u'_{1}=(u_{1}-\sum_{j<l_{1}}a_{1j}x_{r}^{\frac{j}{d}})/
  x_{r}^{\frac{l_{1}}{d}} , \] then
  \(\overline{K'}(u_{1}, x_{r}^{\frac{1}{d}})= \overline{K'}(u'_{1}, x_{r}^{\frac{1}{d}}) \).
  As \( u'_{1} (mod\ {\frak{m}}_{1}) = a_{1j_{1}} \), it follows that  
  \( R_{1}/{\frak{m}}_{1}\supset \overline{K'}(a_{1l_{1}}) \).
  To prove the opposite inclusion, take an element 
  \( h(u'_{1},x_{r}^{\frac{1}{d}})/ g(u'_{1},x_{r}^{\frac{1}{d}})\in 
  R_{1}  \), then the leading coefficients of \( {\alpha'}^*(h) \) 
  and \( {\alpha'}^*(g) \) are in \( \overline{K'}[a_{1l_{1}}] \).
  Therefore the leading coefficient  of   \( {\alpha'}^*(h/g)  \)
  is in   
  \( \overline{K'}(a_{1l_{1}}) \), which shows that the class \( 
  \overline{h/g}\in R_{1}/{\frak{m}}_{1} \) is in \( 
  \overline{K'}(a_{1l_{1}}) \).
  
  Next, assume that \( a_{ij}\in \overline{K'(a_{1l_{1}})} \) for 
  every \( i\geq 2 \), and \( j \).
  Then,
\[    (L[[t']]\cap K(Y) )/((t')\cap K(Y))\subset 
  \overline{K'(a_{1l_{1}})}[[t']]/(t')=\overline{K'(a_{1l_{1}})}. \] 
   Therefore  \[ \trdeg_{\bC}(L[[t']]\cap K(Y)) /((t')\cap K(Y))= n-r+1 < n-1, \] 
  which is a contradiction to Proposition \ref{divisorial}.
   Hence, there exists \( a_{ij} \) \( (i\geq 2) \) transcendental over \( 
   K'(a_{1l_{1}}) \).
   Let it be \( a_{2l_{2}} \).
   By continueing this procedure, we obtain finally 
   \( R_{r-1}=L[[t']]\cap 
   \overline{K'}(u_{1},..,u_{r-1},x_{r}^{\frac{1}{d}})  \) and
   \( R_{r-1}/{\frak{m}}_{r-1}= \overline{K'}(a_{il_{1}}, 
   ..,a_{r-1l_{r-1}}) \), with \( a_{il_{1}}, 
   ..,a_{r-1l_{r-1}} \) algebraically independent over \( K' \).

  (ii)\( \Rightarrow \)(i). Let \( X \) be of dimension \( n \).
   We identify \( K(X) \) with the subfield 
  of \( K((t)) \) by the field homomorphism \( \alpha^* \) induced from 
  \( \alpha \).
  Let \( R=K[[t]]\cap K(X) \) and \( {\frak{m}} \) be the maximal 
  ideal of \( R \).
  Then, by Lemma \ref{divisorial}
  it is sufficient to  prove that \( \trdeg_{\bC}R/{\frak{m}}\geq n-1 \).
  As \( \trdeg_{\bC}K'=n-r \), it is sufficient to prove that 
  \( \trdeg_{K'}R/{\frak{m}}\geq r-1 \).              
  By \( \alpha^*(x_{r})={t'}^d \), there is a \( d \)-th root  \( 
y \) of \( 
  x_{r} \) such that \( \alpha^*: K(X) \to K((t)) \) is lifted as 
  \( \tilde{\alpha}^*: K(X)(y)\to L((t')) \) with \( y \mapsto t' \).
  Let \( \tilde R=L[[t']]\cap K(X)(y)  \) and \( \tilde {\frak{m}} 
  \) the maximal ideal of \( \tilde R \).
  Then, by Lemma \ref{algebraic}, it is sufficient to prove 
  \( \trdeg_{K'}\tilde R/\tilde{\frak{m}}\geq r-1 \).

  By the assumption of (iii) and Lemma \ref{minimal} may assume that  
  \begin{enumerate}
\item[(1)]
\( a_{1j}\in \overline{K'} \) for every \( j< l_{1} \);
\item[(2)]
\( a_{2j}\in \overline{K'(a_{1l_{1}})} \) for every \( j< l_{2} \);
\newline
\ldots
\item[\( (r-1) \)]
\( a_{r-1j}\in \overline{K'(a_{1 l_{1}},\ldots, a_{r-2 l_{r-2}})} \) 
for every \( j<l_{r-1} \),
\end{enumerate}           
  Then identifying an element of \( K(X)(y) \) with the element in \( 
  L((t')) \) by 
  \( \tilde{\alpha}^* \), we have
   \[(\ref{thm}.1)\ \ \ \frac{x_{1}-a_{1k_{1}}y^{k_{1}}-..-a_{1l_{1}-1}y^{l_{1}-1}}{y^{l_{1}}}=
   a_{1l_{1}}+ a_{1l_{1}+1}t' + a_{1l_{1}+2}{t'}^2+\ldots.  \]
  The left hand side is in \( R_{1}:=K(X)(y, a_{1k_{1}},.,a_{1l_{1}-1})
  \cap L[[t']] \).
  Let the maximal ideal of \( R_{1} \) be \( {\frak{m}}_{1} \).
  Then, since 
  \( R_{1}/{\frak{m}}_{1}=\tilde R/\tilde{\frak{m}}(a_{1k_{1}},.,a_{1l_{1}-1}) \), 
  \( R_{1}/{\frak{m}}_{1} \) is 
  algebraic over \( \tilde R/\tilde{\frak{m}} \).   
  The right hand side of (\ref{thm}.1) shows that 
  \( a_{1l_{1}}\in  R_{1}/{\frak{m}}_{1} \).
  
  Next, we obtain      
  \[(\ref{thm}.2)\ \ \ 
  \frac{x_{2}-a_{2k_{2}}y^{k_{2}}-..-a_{2l_{2}-1}y^{l_{2}-1}}{y^{l_{2}}}=
   a_{2l_{2}}+ a_{2l_{2}+1}t' + a_{2l_{2}+2}{t'}^2+\ldots.  \]
 The left hand side is in \( R_{2}:=K(X)(y, a_{1k_{1}},.,a_{1l_{1}-1}, 
 a_{2k_{2}},.,a_{2l_{2}-1})
  \cap L[[t']] \). 
  As \( 
  R_{2}/{\frak{m}}_{2}=R_{1}/{\frak{m}}_{1}(a_{2k_{2}},.,a_{2l_{2}-1}) 
  \),
  this is algebraic over \( R_{1}/{\frak{m}}_{1} \).
  The right hand side of (\ref{thm}.2) shows that 
  \( a_{2l_{2}}\in  R_{2}/{\frak{m}}_{2} \).
  By these successive procedure, we obtain a sequence of algebraic 
  extensions of fields:
  \[ \tilde R/\tilde{\frak{m}}\subset R_{1}/{\frak{m}}_{1} 
  \subset ..\subset R_{r-1}/{\frak{m}}_{r-1}, \]
  with \( a_{1l_{1}},.,a_{r-1l_{r-1}}\in R_{r-1}/{\frak{m}}_{r-1} \).
  Therefore \( \trdeg_{K'}R_{r-1}/{\frak{m}}_{r-1}\geq r-1 \), which 
  implies \( \trdeg_{K'}\tilde R/\tilde{\frak{m}}\geq r-1 \).   

\end{pf}

\begin{cor}
\label{cor}
  Let \( \alpha:\sT\to X \) be a fat arc of a variety \( X \)  and
  \( e \) the center \( \alpha(0) \) of  the arc \( \alpha \)
  with \( \codim _{X}\overline{\{e\}}= r \).
  Let \( \alpha^*: {\o}_{X,e}\to K[[t]] \) be the local homomorphism 
  corresponding to \( \alpha \).
  Assume that there exist \( x_{1},\ldots,x_{r}\in {\frak{m}}_{X,e} \) such 
     that, for \( i=1,\ldots,r-1 \),
     we have \( \alpha^*(x_{i})= \sum_{j}b_{ij}{t}^j\) 
     with  a set of coefficients \( 
     b_{1l_{1}},b_{2l_{2}},\ldots,b_{r-1l_{r-1}} \) 
     which are algebraically 
     independent over \( K'(b_{rj}\ (j<d+\max\{l_{i}\}), b_{ij} (j<{l_{i}}, 
     i=1,.,r-1))
      \),
     where \( K'=\o_{X,e}/{\frak{m}}_{X,e} \) and 
     \( \alpha^*(x_{r})= \sum_{j\geq d}b_{rj}{t}^j\) \( (b_{rd}\neq 0) \).
     Then, \( \alpha \) is divisorial.      
 
\end{cor}

\begin{pf}
  Let  \( L[[t']]\supset K[[t]] \)  be a canonical extension 
  with respect to \( \alpha^*(x_{r})=\sum_{j\geq d}b_{rj}{t}^j \).
  If we write \( t=c_{1}t'+c_{2}{t'}^2+\ldots \), 
  then \( c_{k}\in \overline{\bC(b_{rj}\mid j<d+k)} \),
  where \( \overline{*} \) is the algebraic closure of \( * \) in \( L \).
  If we write \( \alpha^*(x_{i})= \sum_{j}a_{ij}{t'}^j\),
  then \( a_{il_{i}}=b_{il_{i}}c_{1}^{l_{i}}+ d_{i} \),
  where \( d_{i}\in \overline{\bC(b_{rj}\ (j<d+l_{i}), b_{ij}\ (j< 
  l_{i}))} \).
  By the assumption on \( b_{1l_{1}},..,b_{r-1l_{r-1}} \),
  the coefficients 
  \( a_{1l_{1}},..,a_{r-1l_{r-1}} \) are algebraically independent 
  over \( K'(b_{rj}\ (j<d+\max\{l_{i}\}), b_{ij}\ (j< 
  l_{i}, i=1,..,r-1)) \), and therefore algebraically independent over \( K' \).
  Then, we can apply Theorem \ref{thm}.   
\end{pf}

  The following example shows a finite dimensional divisorial arc.

\begin{exmp}
\label{finite}
   Let \( K=\bC(a_{1},..,a_{n}) \), where \( a_{1},..,a_{n} \) are 
   algebraically independent over \( \bC \).
  Let \( \alpha:\sT\to \bA_{\bC}^n=\spec \bC[x_{1},..,x_{n}] \) be an 
  arc defined by \( \alpha^*:\bC[x_{1},..,x_{n}]\to K[[t]] \),
  \( x_{i}\mapsto a_{i}t^{v_{i}} \).
  Then, by Corollary \ref{cor} \( \alpha \) is a divisorial arc.
  We can also see that the corresponding divisorial valuation is a toric 
  valuation  \( \val_{ D_{v}} \),
  where \( v=(v_{1},\ldots,v_{n}) \). 
  To show this, we use the notation and terminologies in \cite{fulton}
  (see also Example \ref{toric}).
  First for a monomial \( x^u \) \( (u=(u_{1},..,u_{n})\in 
  \sigma^{\vee}\cap M )\), it follows \( \ord \alpha^*(x^u)=
  \langle v, u \rangle = \val _{D_{v}}( x^u) \). 
  Here, \( \sigma \) is the positive octant in \( N_{\bR} \)
   defining the toric variety \( 
  \bA_{\bC}^n \).   
  For an element \( \sum_{u}b_{u}x^u\in \bC[x_{1},..,x_{n}] \),
  it is clear that
  \[ (\ref{finite}.1)\ \ \ \ \ \ \ \ \ \ \ \ \ \ \ \ \   
  \ord \alpha^*( \sum_{u}b_{u}x^u)\geq \min_{b_{u}\neq 0}\langle v, 
  u \rangle =\val_{D_{v}}( \sum_{u}b_{u}x^u).\ \ \ \ \ \ \ \ \ \ \ \ \ \ \ \ \ \ \ \]
  If the equality in (\ref{finite}.1) does not hold, then
  \[ \sum_{u} b_{u}(\prod_{i=1}^n a_{i}^{u_{i}})=0, \]
  where the sum is over all \( u \) such that \( b_{u}\neq 0 \) and 
   \( u \) attains \( \min_{b_{u'}\neq 0}\langle v, 
  u' \rangle \).
  This equality gives an algebraic relation of \( a_{1},..,a_{n} \) 
  over \( \bC \), which is a contradiction.  
\end{exmp}

\begin{exmp}[a cylinder on a non-singular variety \cite{ELM}]
\label{cylinder2}
  Let \( X \) be a non-singular variety of dimension \( n \) and \( C 
  \) an irreducible cylinder, i.e., \( C=\psi_{m}^{-1}(S) \) for an 
  irreducible constructible set \( S\subset X_{m} \).
  Here, \( \psi_{m}:X_{\infty}\to X_{m} \) is the morphism of 
  truncation. Then the valuation \( \val_{C} \) defined in \cite{ELM} 
  is divisorial. 
\end{exmp}

\begin{pf}
  Let \( \alpha\in C \) be the generic point, then
  it is sufficient to prove that \( \alpha \) is divisorial. 
  Note that \( \alpha_{m}=\psi_{m}(\alpha) \) is the generic point of \( S \).
  Let the codimension of the center \( e \) of \( \alpha \) be \( r \). 
  Then, by taking a suitable open neighborhood of the center, 
  we may assume that the closure of the center is defined by \( r 
  \)-functions \( x_{1},\ldots,x_{r} \).  
  Let \( \alpha_{m}^*:\o_{X,e} \to K_{m}[t]/(t^{m+1}) \) be the ring 
  homomorphism corresponding to  \( \alpha_{m}\). 
  Then \( \alpha_{m}^* \) can be extended to a local homomorphism 
  \( \alpha_{m}^*:\widehat{\o_{X,e}}=K'[[x_{1},\ldots,x_{r}]]\to  
  K_{m}[t]/(t^{m+1})  \),
  where \( K' \) is the residue field of \( \o_{X,e} \).
  Therefore, \( \alpha_{m}^* \) is determined by \( K'\hookrightarrow 
  K_{m} \) and the images of \( x_{i} \)'s.  
  If \( \alpha_{m}^*(x_{i})=\sum _{j=1}^ma_{ij}t^j \) \( (a_{ij}\in 
  K_{m}) \),
  \( \alpha \) is given by 
   \( \alpha^*:\widehat{\o_{X,e}}\to K[[t]] \)
  with the inclusion \( K'\hookrightarrow K_{m}
  \hookrightarrow K \) and the images    
  \( \alpha^*(x_{i})=\sum_{j=1}^{\infty}a_{ij}t^j .\)
   Here, \[ K=K_{m}(a_{i,m+1},a_{i,m+2},\ldots\mid i=1,\ldots,r) \] and 
   \( a_{i,m+1},a_{i,m+2},\ldots \) \( (i=1,\ldots,r) \) are 
   algebraically independent over \( K_{m} \). 
  Hence, by Corollary \ref{cor} \( \alpha \) is divisorial.   
\end{pf}

\begin{exmp}[divisorial  sets 
 in 
\cite{d-l}, McKay correspondence in \cite{i-r}]
  Let \( d  \geq 1 \) be an integer and \( G \) a finite subgroup of \( 
  \GL_{n}(\bC) \) of order \( d \).
  We fix a primitive \( d \)-th root of unity \( \zeta\in \bC \).
  Let \( X \) be the quotient of \( \bA_{\bC}^n \) by the action of \( 
  G \) and \( h:\bA_{\bC}^n\to X \) the canonical projection.
  The following construction of a subset \( (X_{\infty}^0)[g] \) of   \( 
  X_{\infty} \) corresponding to a conjugacy class \( [g] \) of \( G \) 
  is due to
 J. Denef and F. Loeser (\cite[2.1]{d-l}).
  We denote the origin of \( \bA_{\bC}^n \) and its image in \( X \) 
  by \( 0 \).
  Let \( X_{\infty}^0 =(X_{\infty}\setminus (\sing X)_{\infty})\cap
  \pi^{-1}_{X}(0)\).
  Let \( \alpha\in X_{\infty}^0 \) be the arc
  \( \sT \to X \) and \( \overline{K} \) the algebraic closure of \( K \).
  We denote the induced arc \( \spec \overline{K}[[t]]\to X \) 
   again by \( \alpha \).
  Then  we can lift \( \alpha \) to 
  a morphism \( \tilde\alpha \) making the following diagram 
  commutative:
\[ \begin{array}{ccc}
\Spec \overline{K}[[t^{1/d}]]& \stackrel{\tilde\alpha}\longrightarrow& 
\bA_{\bC}^n\\
\downarrow & & \downarrow h\\
\spec \overline{K}[[t]]& \stackrel{\alpha}\longrightarrow& X.\\
\end{array} \] 
  For an element \( g\in G \) let 
   \[ (X_{\infty}^0)_{g}=\{\alpha\in X_{\infty}^0 \mid 
   \tilde{\alpha}(\zeta t^{1/d})=g\tilde{\alpha}(t^{1/d})\}. \]
  Then \( g \) and \( g' \) are conjugate if and only if  
  \( (X_{\infty}^0)_{g}=(X_{\infty}^0)_{g'} \).
  Hence, we can define the subset \( (X_{\infty}^0)[g]:=(X_{\infty}^0)_{g} \)
  for a conjugacy class \( [g] \) in \( G \).
  We have a decomposition
  \[ X_{\infty}^0=\coprod_{[g]} (X_{\infty}^0)[g], \]
  with   each \(  (X_{\infty}^0)[g] \)   irreducible.  

  We are going to show that \( (X_{\infty}^0)[g] \) is a divisorial set.  
  For \( g\in G \) taking a suitable coordinates system \( 
  x_{1},\ldots,x_{n} \) of \( 
  \bA_{\bC}^n \),
  we may assume that the matrix  \(g
   \) is diagonal and the \( i \)-th diagonal coefficient is \( 
   \zeta^{e_{i}} \) with \( 1\leq e_{i}\leq d \)  \( (i=1,\ldots,n) \).
  Then, we have a homomorphism 
  \[\Lambda: \bC[x_{1},\ldots,x_{n}]
  \to A_{\infty}[[t^{1/d}]] ,\ \ 
   x_{i}\mapsto t^{e_{i}/d}\alpha^*(x_{i})   ,\] 
  where \( \spec A_{\infty}=(\bA_{\bC}^n)_{\infty} \) and 
  \( \alpha^*:\bC[x_{1},\ldots,x_{n}]\to A_{\infty}[[t]] \) is the 
  ring homomorphism corresponding to the universal arc on \( 
  \bA_{\bC}^n \).
  Here, we note that  
  \[ \alpha^*(x_{i})=\sum_{j=0}^{\infty}a_{ij}t^j, \]
  where \( \{a_{ij}\}_{1\leq i \leq n, j\geq 0} \) 
  are algebraically independent over \( \bC \).  
  By restricting \( \Lambda \) to the subring, we obtain a homomorphism 
  \( {\lambda'}^*:\bC[x_{1},\ldots,x_{n}]^{\langle g \rangle } \to  A_{\infty}[[t ]] 
  \) whose restriction also gives  a homomorphism 
  \[ \lambda^*: \bC[x_{1},\ldots,x_{n}]^G \to  A_{\infty}[[t ]] 
  .\]
  Let \( K \) be the quotient field of \( A_{\infty} \).
  Then, the center of \( \lambda:\sT\to X \) is \( 0 \) and \( \lambda 
  \) factors through the generic point \( \gamma: \spec 
  \bC(\gamma)[[t]]\to X \) of \( (X_{\infty}^0)[g] \) (\cite[2.3.4]{d-l}).
   Therefore, in order to show that \( \gamma \) is divisorial 
  it is sufficient to prove that \( \lambda \) is divisorial.  
  By Lemma \ref{commutative}, it is also sufficient to show that 
  \( \lambda':\sT\to \bA_{\bC}^n/\langle g \rangle \) is divisorial.
  The center of \( \lambda' \) is \( 0 \) in   \( \bA_{\bC}^n/\langle g 
  \rangle \) and, for \( i=1,.., n \), \( x_{i}^d \in 
  \o_{\bA_{\bC}^n/\langle g 
  \rangle, 0} \) is mapped to \( 
  t^{e_{i}}(a_{i0}+a_{i1}t+a_{i2}t^2+..)^d \) by \( {\lambda'}^* \),
  where \( K=\bC(a_{i0},a_{i1},a_{i2},..\mid i=1,..,n) \) and \( a_{ij} 
  \)'s are algebraically independent over \( \bC \).
  By Corollary \ref{cor} \( \lambda' \) is divisorial.
  
  To see the concrete correspondence 
  \( (X_{\infty}^0)[g]\mapsto v_{\gamma} \),
  let \( N=\bZ^n \) be the lattice for a toric variety \( \bA_{\bC}^n \). 
  Then \( N'=N+\frac{1}{d}(e_{1},..,e_{n})\bZ \) is the 
  lattice for a toric variety \( \bA_{\bC}^n/\langle g \rangle  \).
  If we put \( v=\frac{1}{d}(e_{1},..,e_{n})\in N' \),
  it follows that \[\ord_{t}{\lambda'}^*(f) =\min _{x^u\in f}  
  \langle v, u\rangle  
  \]for a regular function \( f \) on \( \bA_{\bC}^n/\langle g \rangle  \).
  The proof is the same as in Example \ref{finite}.
  Therefore the valuation \( v_{\lambda'} \) is \( \val_{D_{v}} \).
  Hence, the valuation \( v_{\gamma}=v_{\lambda} \) is the restriction of 
  the toric divisorial valuation \( \val_{D_{v}} \).
  Consider the case \( G\subset \SL_{n}(\bC) \) and \( \dim X=3 \).
  Restricting the map   \( [g]\mapsto (X_{\infty}^0)[g] \mapsto 
  \val_{D_{v}} \) onto the subset \( \Gamma_{1}^0 \) consisting of 
  \( [g] \)'s with \( \sum e_{i}=d \),
  we obtain the ``McKay correspondence'' in \cite[Theorem 1.6]{i-r}. 
\end{exmp}

\begin{cor}
  Let \( G \) be a finite abelian subgroup of \( \GL_{n}(\bC) \) acting 
  on \( \bA_{\bC}^n \).
  Assume that \( X=\bA_{\bC}^n/G  \) has an isolated singularity.
  Then, 
  \begin{center}
  \( \# \{ \) essential divisors over \( X \}\leq 
  \#G-1\) 
  \end{center}
\end{cor}

\begin{pf} 
  Here we also fix a primitive \( d \)-th root of unity \( \zeta \),
  where \( d=\# G \).
  As \( G \) is abelian, every conjugacy class consists of only one 
  element of \( G \).
  Every element \( g\in G \) can be written as a diagonal matrix 
  with the \( i \)-th diagonal coefficient \( \zeta^{e_{i}} \)
  with \( 1\leq e_{i}\leq d \).
  We write \( v_{g}=\frac{1}{d}(e_{1},\ldots,e_{n})\in N' \),
  where \( N' \) is the lattice for a toric variety \( X \).
  First we prove that 
  \[ (X_{\infty}^0)[1]\subset \overline{(X_{\infty}^0)[g]} \] 
  for some \( g\neq 1 \in G \).
  Since
   \[ X_{\infty}^0=\coprod_{[g]} (X_{\infty}^0)[g], \]
   the irreducible components of \(  X_{\infty}^0 \) are the closures 
   of \( (X_{\infty}^0)[g] \)'s for some \( g\in G \).
  If there is an irreducible component \( \overline{(X_{\infty}^0)[g]} 
  \) with \( g\neq 1 \),
  then its generic point is the generic point of \( 
  \varphi_{\infty}(\pi_{Y}^{-1}(D_{v_{g}}) )\) (c.f., Example 
  \ref{nash}), therefore is 
  the generic point of \( T_{\infty}(v_{g}) \) by Lemma \ref{equal}.
  The generic point of  \( (X_{\infty}^0)[1]  \) belongs to \( 
  T_{\infty}(v_{1}) \) which is contained in the closure of 
  \( T_{\infty}(v_{g}) \), because \( v_{1}\geq v_{g} \)
  (\cite[Proposition 4.8]{i}, see  Proposition \ref{dominate} in 
  this paper).
  If \( \overline{(X_{\infty}^0)[1]} \) is only irreducible component 
  of \( X_{\infty}^0 \), then it must be the closure of 
  \( T_{\infty}(v_{1})  \).
  Hence, \( \overline{T_{\infty}(v_{1})}\supset T_{\infty}(v_{g})  \) 
  for \( g\neq 1 \)
  which is a contradiction to \( v_{1}\geq v_{g} \) 
  (\cite[Proposition 4.8]{i}).
  
  Now we obtain that the number of irreducible components of 
  \( X_{\infty}^0 \) is less than \( d \).
  As \( X \) has an isolated singularity at \( 0 \) the components of 
  \( X_{\infty}^0 \) are the Nash components. 
  As Nash map is bijective for a toric variety(\cite{i-k}), 
  the number of the essential divisors is less than \( d \).
\end{pf}
  T. Mizutani \cite{mizutani} proved this corollary by an elementary 
  way and gave  examples
  that there are  
  exactly \( d-1 \) essential divisors.   
\vskip2truecm

\section{Essential divisors and the arc space of a toric variety}

\noindent
  In this section we summarize the notion of essential divisors and
  basic properties of the arc space of a toric variety.

  When we treat a toric variety,
  we use the terminologies in \ref{toric}.

\begin{defn} Let $X$ be a variety, \( \psi:X_1 \to X \)
a proper birational morphism from a normal variety \( X_{1} \)
 and $E\subset X_1$ 
an irreducible exceptional divisor  of 
  \( \psi \). 
 Let  \( \varphi:X_2 \to X \) be another
 proper birational morphism from a normal variety \( X_{2} \).
The birational  map
  \(\varphi^{-1}\circ \psi  : X_1  \dasharrow  X_2 \) is 
defined on a (nonempty) open subset $E^0$ of $E$.
The closure of $(\varphi^{-1}\circ \psi)(E^0)$ is well defined.
It is called the {\it center} of $E$ on $X_2$.

  We say that \( E \) appears in \( \varphi \) (or in \( X_2 \)),
  if 
the center of $E$ on $X_2$ is also a divisor. In this case
the birational  map
  \(\varphi^{-1}\circ \psi  : X_1  \dasharrow  X_2 \) is  a local isomorphism at the 
  generic point of \( E \) and  we denote the birational
  transform of \( E \) on \( X_2 \)
  again by \( E \). For our purposes $E\subset X_1$ is identified
with $E\subset X_2$. 
 Such an equivalence class  is called an {\it exceptional divisor over $X$}. 
\end{defn}

\begin{defn}
  Let \( X \) be a variety over \( \bC \). In this paper, by  a
{\it  resolution} of the 
  singularities of \( X \) 
   we mean  a proper, birational  morphism \( \varphi:Y\to X \) 
with  \( Y \)  non-singular
such that $Y\setminus \varphi^{-1}(\sing X)\to X\setminus \sing X$
 is an isomorphism.
 Here, \( \sing X \) is the singular locus of \( X \).

 \end{defn}
  
  \begin{defn}\label{nashessential}
  An exceptional divisor \( E \) over \( X \) is called an {\it essential 
  divisor} over \( X \)
   if for every resolution \( \varphi:Y\to X \)
the center of \( E \) on \( Y \) is an irreducible component of 
   \( \varphi^{-1}(\sing X) \).
   
   For a toric variety \( X \), an {\it equivariant essential divisor}
   over \( X \) is a divisor \( E \) over \( X \) whose center on every 
   equivariant resolution \( \varphi:Y\to X \)  is an irreducible component of 
   \( \varphi^{-1}(\sing X) \).

\end{defn}

About  essential divisors, we have the following:

\begin{prop}[\cite{i-k}]
\label{charac-essential}
  For a toric variety \( X \) the notions an essential divisor and an 
  equivariant essential divisor coincide.
  
  Let \( \sigma \) be the cone defining an affine toric variety \( X 
  \) and \( v\in \sigma\cap N \).
  Then \( D_{v} \) is essential if and only if 
  \( v \) is a minimal element in \[ S=\cup _{\tau< \sigma 
  :{\operatorname{singular}}}\tau^o \cap N ,\]
   with respect to the order \( 
  \leq_{\sigma} \) ( \( v\leq_{\sigma}v' \Leftrightarrow v'-v\in \sigma 
  \) ). 
  Here, \( \tau^o \) means the relative interior of \( \tau \).
\end{prop}

  In \cite{i}, we introduce a locally closed subset \( T_{\infty}(v) \)
   of the arc space 
  \( X_{\infty} \) of an affine toric variety \( X \) as follows (see 
  \ref{toric}):
  \[  T_{\infty}(v)= \{\alpha\in X_{\infty} \mid \alpha(\eta)\in T,\  \ord_{t} 
  \alpha^*(x^u)=\langle v, u \rangle \ {\operatorname{for}}\ u\in M \}. \]
  In order to exhibit the space \( X \) explicitly, we 
  denote \( T_{\infty}(v) \) by \( T_{\infty}^X(v) \).
  The following is obtained in \cite{i}.
   
\begin{prop}[\cite{i}]
\label{dominate}
  Let \( X \) be an affine toric variety defined by a cone \( \sigma \) 
  in \( N \).
  For \( v,v' \in \sigma\cap N  \),
  the relation \( v\leq_{\sigma} v' \) holds if and only if   
    \( \overline{T_{\infty}^X(v)}\supset T_{\infty}^X(v') \).

\end{prop}

\begin{lem}
\label{dominate2}
  Let \( X \) be an affine toric variety defined by a cone \( 
  \sigma  \) in \( N \) and  \( v, v'\) elements in \( 
  \sigma\cap N \).
  Let \( \varphi:Y \to X \) be an 
  equivariant proper birational morphism in which \( D_{v} \) and \( D_{v'} \) 
  appear.
  If \( v\leq_{\sigma}v' \), then \( {\varphi(D_{v})}\supset
  \varphi(D_{v'}) \).
\end{lem}

\begin{pf}
  If \( v\in \tau^o \), \( v'\in {\tau'}^o \) for faces \( 
  \tau,\tau'<\sigma \), then 
  \( \varphi(D_{v})=\overline{orb \tau} \) and 
  \(  \varphi(D_{v'})=\overline{orb \tau'} \).
  By the assumption of the lemma, it follows that \( v'=v+v'' \) for 
  some \( v''\in \sigma \).
  As \( v'\in \tau' \), \( v, v''\in \tau' \).
  Hence, \( \tau<\tau' \),  which yields the assertion of the lemma.  
\end{pf}

At the end of this section, we prove a technical lemma which is used 
in the next section.

\begin{lem}
\label{miz}
  Let \( \sigma=\langle e_{1},\ldots,e_{m}\rangle \) be a singular simplicial 
  cone in \( N \).
  Assume that one facet is non-singular, then \( v=\sum_{i=1}^mb_{i}e_{i} \) 
  \( (b_{i}\geq 1,\ i=1,..,m) \) is 
  not minimal in \( S=\cup _{\tau< \sigma 
  :{\operatorname{singular}}}\tau^o \cap N \).
\end{lem}

\begin{pf}
  We may assume that \( e_{1}=(1,0,..,0), e_{2}=(0,1,0,.,0),..,
  e_{m-1}=(0,.,0,1,0),  \) and \( e_{m}=(a_{1}c-c_{1}, a_{2}c-c_{2},..,
  a_{m-1}c-c_{m-1},c) \), where \( a_{i} \)'s are integers and 
  \( 0\leq c_{i}< c \).
  Here, we note that \( c>1 \), since \( \sigma \) is singular.
  Then, \( 
  v'=\frac{c_{1}}{c}e_{1}+\frac{c_{2}}{c}e_{2}+..+\frac{c_{m-1}}{c}e_{m-1}+
  \frac{1}{c}e_{m} 
  \) is in \( S \) and \( v'\leq_{\sigma} v \).  
\end{pf}

\section{the Nash problem for a pretoric variety}

\begin{defn}
\label{pretoric}
  A variety \( X \) is called a {\it pretoric variety} if 
  \begin{enumerate}
  \item[(1)]
  there are a toric variety \( Z \) with the torus \( T' \) and 
  a finite morphism \( \rho:X\to Z \) \'etale on \( T' \),
  \item[(2)]
  for the normalization \( \nu:\overline{X} \to X \), 
  \( \overline{X} \) is a toric variety with the torus \( T \) and 
  the composite \( \rho\circ \nu:\overline{X}\to Z \) is the 
  equivariant quotient morphism by the group \( N'/N \), where \( N 
  \) and \( N' \) are the lattice on which the fans of \( \overline{X} 
  \) and \( Z \) are defined, respectively, and 
  \item[(3)]
  the subset \( \nu^{-1}(\sing X) \) is an invariant closed set on
  \( \overline{X} \).
\end{enumerate}  
\end{defn}

  We will see two typical examples of a pretoric variety.

\begin{say}[\cite{gkz}]
  When we say ``a toric variety'', it means always a normal toric 
  variety. 
  Here, we introduce a not-necessarily normal affine toric variety.
  A not-necessarily normal affine toric variety is of the form
   \( X_{ \Gamma }= \spec \bC[ \Gamma ] \), where \( \Gamma \subset M= \bZ^n\) is a finitely 
  generated semigroup with \( 0 \) and \( \Gamma \) generates the abelian 
  group \( M \).
  Then, the torus \( T=\spec \bC[M] \) acts on \( X_{ \Gamma } \).
  Denote by \( K( \Gamma )\subset  M_{\bR}  \), the convex cone which is the 
  convex hull of \( \Gamma \) and by \( \overline{ \Gamma } \) the intersection 
  \( K(\Gamma)\cap M \).
  Then, \( X_{\overline{ \Gamma }} \) is a normal toric variety and 
   the the inclusion \( \bC[ \Gamma ]\hookrightarrow \bC[\overline{ \Gamma }] \)
  induces the equivariant normalization \( X_{\overline{ \Gamma }}\to X_{ \Gamma } \).
\end{say}

\begin{exmp}
  A not-necessarily normal toric variety is a pretoric variety.
  This is proved as follows:
  Let \( X=\spec \bC[\Gamma] \) be a not-necessarily normal toric 
  variety of dimension \( n \) and \( \overline{X}=\spec \bC[\sigma^{\vee}\cap M] \) the 
  normalization of \( X \).
  Subdivide \( \sigma^{\vee} \) into simplicial cones without adding 
  any 1-dimensional cones.
  Let \( \tau_{1}, \tau_{2},.., \tau_{s} \) be the \( n 
  \)-dimensional simplicial cones which are obtained by this subdivision.
  We can take  generators \( e_{1}^{(i)},..,e_{n}^{(i)}
  \) of \( \tau_{i} \) in \( \Gamma \).
  Define \( M_{i}=\oplus_{j=1}^n\bZ e_{j}^{(i)} \), then 
  \(  M_{i} \) is a subgroup of \( M \) of finite index.
  Let \( M' \) be the intersection \( \bigcap_{i=1}^s M_{i} \).
  Then, \( M' \) is a subgroup of \( M \) of 
  finite index.
  It follows that \( \sigma^{\vee}\cap M' \subset \Gamma \).
  Indeed, an arbitrary element \( u\in \sigma^{\vee}\cap M' \) is 
  contained in \( \tau_{i}\cap M_{i} \) for some \( i \).
  Then, by the definition of \( M_{i} \), we have 
  that \( u=\sum_{j=1}^n a_{j}e_{j}^{(i)} \) with \( a_{j}\in 
  \bZ_{\geq 0} \).
  As \( e_{j}^{(i)} \)'s are in \( \Gamma \), it follows that \( u\in 
  \Gamma \).
  By this inclusion  \( \sigma^{\vee}\cap M' \subset \Gamma  \) we 
  obtain a finite morphism \( \rho: X\to Z=\spec 
  \bC[\sigma^{\vee}\cap M'] \).
  The other conditions for a pretoric variety follows immediately.    
\end{exmp}

The following is an example of a pretoric variety without a toric 
action.

\begin{exmp}
 Let \( \overline{X} \) be \( \spec \bC[x,y] \) and \( X \) be 
 \( \spec \bC[x,y^3, y^4] \), then \( X \) is a non-normal toric 
 variety with the normalization \( \nu:\overline{X}\to X \).
 Therefore we have a diagram \( 
 \overline{X}\stackrel{\nu}\longrightarrow X 
 \stackrel{\rho}\longrightarrow Z \) as in Definition \ref{pretoric}.
 Let \( X_{0} \) be \( \spec \bC[x, y+y^2, y^3, y^4] \),
 then \( X_{0} \) is a pretoric variety with the diagram:
 \( \overline{X}\to X_{0}\to Z \).
 By the definition, \( X_{0} \) does not admit a toric action.
\end{exmp}

\begin{defn}
\label{pretorict}
  Let \( X \) be an affine pretoric variety and 
  \( \overline{X}\stackrel{\nu}\longrightarrow X 
  \stackrel{\rho}\longrightarrow Z \) the diagram as in the definition 
  \ref{pretoric}.
  Let \( T \) and \( T' \) be the tori of \( \overline{X} \) and 
  \( Z \), respectively.
  Let \( \overline{X} \) and \( Z \) be defined by a cone \( \sigma \) 
  in \( N \) and \( \sigma \) in \( N' \), respectively.
  The dual of \( N, N' \) is denoted by \( M, M' \). 
  For \( v\in \sigma\cap N \subset \sigma \cap N' \) the subsets 
  \( T_{\infty}^{\overline{X}}(v) \) and \( T_{\infty}^{{Z}}(v)  \) 
  are defined, since \( \overline{X} \) and \( Z \) are toric 
  varieties. 
  Here, we define \( T_{\infty}^{{X}}(v)  \) for a pretoric 
  variety \( X \) as follows:
  \[  T_{\infty}^X(v)= \{\alpha\in X_{\infty} \mid 
  \alpha(\eta)\in\rho^{-1}( T')=\nu(T),\  \ord_{t} 
  \alpha^*(x^u)=\langle v, u \rangle \ {\operatorname{for}}\ u\in M \}, \]
\end{defn}

\begin{lem}
\label{bij}
  Let \( X \) be an affine pretoric variety. 
  Under the same notation as in the previous definition, for \( v\in 
  \sigma\cap N \)
\begin{enumerate}
\item[(i)]
  the morphism \( \nu_{\infty} \) gives a bijection 
  \( T_{\infty}^{\overline{X}}(v)\to T_{\infty}^X(v) \) and
\item[(ii)]
  the morphism \( (\rho\circ\nu)_{\infty} \) gives a surjection 
  \( T_{\infty}^{\overline{X}}(v)\to T_{\infty}^Z(v) \).
\end{enumerate}      
\end{lem}

\begin{pf}
  By the definition of \( T_{\infty}(v) \)'s it is clear that the images 
  of these morphisms are in the target sets.
  For the proof of (i), take an arc \( \alpha\in T_{\infty}^X(v)  \),
  then \( \alpha(\eta)\in \nu(T)\simeq T \) and therefore \( \alpha 
  \) is lifted uniquely to \( \overline{X} \), as \( \nu \) is proper.
  For the proof of (ii), take an arc \( \alpha\in T_{\infty}^Z(v) \).
  Then \( \alpha \) corresponds to a ring homomorphism
  \( \alpha^*: \bC[\sigma^{\vee}\cap M']\to K[[t]] \).
  Let \( \overline{K} \) be the algebraic closure of \( K \) and 
  denote the composite of \( \alpha^* \) and the canonical inclusion
  \( K[[t]]\hookrightarrow \overline{K}[[t]] \) 
  by again \( \alpha^* \).
  Then, by Denef and Loeser \cite{d-l}, we obtain the following 
  commutative diagram for \( d=\# N'/N \):
  \[ \begin{array}{ccc}
  \bC[\sigma^{\vee}\cap M] & \stackrel{\tilde \alpha^*}\longrightarrow &
  \overline{K}[[t^{1/d}]]\\
  \bigcup & & \bigcup\\
   \bC[\sigma^{\vee}\cap M'] & \stackrel{\alpha^*}\longrightarrow &
  \overline{K}[[t]].\\
  \end{array}\]
  As \( v\in \sigma\cap N\subset \sigma\cap N' \), it follows that  \( \langle v, u 
  \rangle \) is a non-negative integer for \( u\in 
  \sigma^{\vee}\cap M \). 
  For \( u\in 
  \sigma^{\vee}\cap M \),
  let \( \tilde\alpha^*(x^u)=\sum_{i=0}^\infty a_{i}t^{i/d}\).   
  By considering \( \alpha^*(x^{du})=\tilde \alpha^*(x^u)^d \), 
  we obtain that  \( a_{i}=0 \) for \( i\not\equiv 0 (mod\ d) \),
  which means \( 
  \tilde{\alpha}^*(x^u)\in \overline{K}[[t]] \).
  This gives the surjectivity of the morphism in (ii).      
\end{pf}

\begin{lem}
\label{pretoricdominate}
   Let \( X \) be an affine pretoric variety. 
  Under the same notation as in Definition \ref{pretorict}, let \( v, v'\in 
  \sigma\cap N \).
  Then, 
  \( v\leq_{\sigma}v' \) if and only if 
  \( \overline{T_{\infty}^X(v)}\supset T_{\infty}^X(v') \).
\end{lem}

\begin{pf}
  If \( v\leq_{\sigma}v' \), then by Proposition \ref{dominate} it 
  follows that 
  \( \overline{T_{\infty}^{\overline{X}}(v)}\supset T_{\infty}^
  {\overline{X}}(v')   \).
  Hence, by Lemma \ref{bij}, (i), we obtain 
  \( \overline{T_{\infty}^X(v)}\supset T_{\infty}^X(v') \).
  Conversely, if \(  \overline{T_{\infty}^X(v)}\supset T_{\infty}^X(v')  \),
  then by Lemma \ref{bij}, (ii), we have that 
  \( \overline{T_{\infty}^Z(v)}\supset T_{\infty}^Z(v') \). 
  Hence, by Proposition \ref{dominate}, it follows that  
  \( v\leq_{\sigma}v' \).
\end{pf}

\begin{say}
\label{notation}
  Let \( X \) be an affine pretoric variety and 
  \( \overline{X}\stackrel{\nu}\longrightarrow X 
  \stackrel{\rho}\longrightarrow Z \) the diagram as in the definition 
  \ref{pretoric}.
  Let \( T \) and \( T' \) be the tori of \( \overline{X} \) and 
  \( Z \), respectively.
  Let \( \overline{X} \) and \( Z \) be defined by a cone \( \sigma \) 
  in \( N \) and \( \sigma \) in \( N' \), respectively.
  The dual of \( N, N' \) is denoted by \( M, M' \). 
  
  As \( \nu^{-1}(\sing X)\) is  an invariant closed set, 
  an irreducible component of 
  \( \overline{\nu^{-1}(\sing X)\setminus \sing \overline{X} } \) is 
  written by \( \overline{ orb \tau_{i}} \)  for some  
  non-singular face  
  \( \tau_{i}<\sigma \) \( (i=1,..,r) \).
  Let \( e_{i} \) be the barycenter of \( \tau_{i} \), 
  i.e., \( e_{i} \) is the sum of the generators of \( \tau_{i} \) in \( 
  N \). 
  
  Let \( v_{j} \) \( (j=1,..,s) \) be the minimal elements of 
  \( S=\bigcup_{\tau<\sigma:
                 \operatorname{\ singular}}   \tau^o \cap N \), i.e., \( D_{v_{j}} \) 
                  \( (j=1,..,s) \) are the essential divisors over \( 
                  \overline{X} \) (see Proposition 
                  \ref{charac-essential}).
  We consider the minimal elements of \( \{e_{i}, 
  v_{j}\}^{i=1,.,r}_{{j=1,.,s}} \) 
  and obtain the following:
\end{say}  
  
\begin{lem}
  Each \( e_{i}  \) \(( i=1,..,r )\) is minimal among \( \{e_{i}, 
  v_{j}\}^{i=1,.,r}_{{j=1,.,s}} \) with respect to the order \( 
  \leq_{\sigma} \). 
\end{lem}

\begin{pf}
  Let \( \varphi:Y\to \overline{X} \) be an equivariant resolution 
  in which \( D_{e_{i}} \)'s and \( D_{v_{j}} \)'s appear.
  Then, \( \varphi(D_{e_{i}})=\overline{orb \tau_{i}} \)  is an 
  irreducible component of  
  \( \overline{\nu^{-1}(\sing X)\setminus \sing \overline{X} } \).
  Therefore, \( \varphi(D_{e_{i}})\not\subset \varphi(D_{e_{k}}) \) 
  for  \( k\neq i \). 
  On the other hand, as \( \varphi(D_{v_{j}})\subset \sing \overline{X} \),
  it follows that \( \varphi(D_{e_{i}})\not\subset \varphi(D_{v_{j}}) \).
  Hence, by Lemma \ref{dominate2}, we obtain that \( 
  e_{k} \not\leq_{\sigma}  e_{i} \) \( (k\neq i) \) and 
  \(v_{j} \not\leq_{\sigma}  e_{i} \) \( (j=1,..,s) \). 
  Thus, \( e_{i} \) is minimal.   
\end{pf} 

\begin{lem}
\label{mainlemma}
  Let \( \{e_{i}, 
  v_{j}\}^{i=1,.,r}_{{j=1,.,w}} \) \( (w\leq s) \) be the set of 
  minimal elements in \( \{e_{i}, 
  v_{j}\}^{i=1,.,r}_{{j=1,.,s}} \).
  Then, there is the inclusion
\[ \{\operatorname{ essential\ divisors\ over\ }{X}\}
\subset  \{D_{e_{i}}, 
   D_{v_{j}}\}^{i=1,.,r}_{{j=1,.,w}}. \] 
\end{lem}

\begin{pf}
  Let \( D \) be an essential divisor over \( X \).
  Let \( \varphi:Y\to X \) be a resolution on which \( D \) appears, 
  then \( \varphi \) 
  factors through the normalization: 
  \[ Y\stackrel{\psi}\longrightarrow \overline{X} \stackrel{\nu}
  \longrightarrow X. \]
  Here, we may assume that \( \psi \) is an equivariant morphism.
  Then, we can put \( D=D_{v} \) for some \( v\in \sigma\cap N \).
  As \( \varphi(D_{v})\subset \sing X \), we have 
  \( \psi(D_{v})\subset \sing \overline{X} \) or 
  \( \psi(D_{v})\subset \overline{orb \tau_{i}} \).

{\bf Case 1.}
\( \psi(D_{v})\not\subset \sing \overline{X} \) and \( \psi(D_{v})\subset 
\overline{orb \tau_{i}} \).

  In this case we show that \( \psi(D_{v})=
  \overline{orb \tau_{i}} \) and \( v=e_{i} \).
  Let \( \psi(D_{v})=\overline{orb \gamma} \), then 
  \( \gamma \) is a non-singular face of \( \sigma \) such that 
  \( \tau_{i}<\gamma \).
  Let \( \psi_{0}:Y_{0}\to \overline{X} \) be an equivariant resolution of 
  \( \overline{X} \).
  As \( \psi_{0 } \) is isomorphic away from \( \sing \overline{X} \),
  \( orb \gamma, orb \tau_{i}  \) are still on \( Y_{0} \).
  Then, take the blow up \( \psi_{1}:Y_{1}\to Y_{0} \) with the center 
  \( \overline{orb \tau_{i}} \), then \( Y_{1} \) is again non-singular.
  Here, \( \psi_{1} \) corresponds to the star-shaped subdivision \( \Sigma_{1} \)
   of the fan \( \Sigma_{0} \)
  of \( Y_{0} \) by \( e_{i} \).
  Take a cone \( \gamma' \) in \( 
  \Sigma_{1} \) such that \( v\in {\gamma'}^o \).
  Here, we note that \( \langle e_{i} \rangle <\gamma' \),
  because, \( v \) is in the relative interior of \( \gamma \)
  and the subdivision is star shaped with the center \( e_{i}\in \gamma \).
  Thus, the center of \( D_{v} \) on \( Y_{1} \) is 
  \( \overline{orb \gamma'} \) which is contained in \( D_{e_{i}} \).
  Therefore, the center of \( D_{v} \) can be a component of 
  \( (\nu\circ\psi_{0}\circ\psi_{1})^{-1}(\sing X) \)  only if  
  \( \overline{orb \gamma'} = D_{e_{i}} \), in which case 
  \( \gamma'=\langle e_{i} \rangle \). 
  This implies \( v=e_{i} \).  
  
  {\bf Case 2.}
  \( \psi(D_{v})\subset\sing \overline{X} \).
  
  In this case, we prove that \( v=v_{j} \) \( (j=1,..,w) \).
  First, \( D_{v} \) must be an essential divisor over \( \overline{X} \).
  Because, if \( D_{v} \) is not an essential divisor over \( 
  \overline{X} \), there is a resolution \( \psi': Y'\to \overline{X} \)
  such that the center of \( D_{v} \) on \( Y' \) is not an irreducible component 
  of \( {\psi'}^{-1}(\sing \overline{X}) \). 
  Therefore the center of \( D_{v} \) cannot be an irreducible component of 
  \( (\nu\circ\psi')^{-1}(\sing X) \).
  Now, for the lemma, it is sufficient to prove that \( D_{v_{j}} \) 
  \( (j>w) \) is not an essential divisor over \( X \).
  For this, we construct an equivariant birational morphism 
  \( \psi'':Y''\to \overline{X} \) such that \( \nu\circ\psi'' \) is 
  a resolution of \( X \) and the center of \( D_{v_{j}} \) on \( Y'' 
  \) is not an irreducible component of   \( 
  (\nu\circ\psi'')^{-1}(\sing X) \).  
  
  As \( v_{j} \) \( (j>w) \) is not minimal in  
  \( \{e_{i}, 
  v_{j}\}^{i=1,.,r}_{{j=1,.,s}} \),
  there is an element \( e_{i} \) such that \( e_{i}\leq_{\sigma}v_{j} \).
  Then, there is an element \( e'_{i}\in \sigma\cap N \) such that 
  \( v_{j}=e_{i}+e'_{i} \).
  Here, \( e'_{i} \) is in the relative interior of a non-singular face \( 
  \tau \) of \( \sigma \).
  This is proved as follows:     As \( v_{j} \in S \)
  (see, \ref{notation}), \( v_{j}\in \gamma^o \) 
  for some singular face \( \gamma<\sigma \).
  If \( e'_{i}\) is in the relative interior of \( \gamma\)  or 
  \( e'_{i} \) is in the relative interior of a singular face of \( 
  \gamma \), then it contradicts to the minimality 
  of \( v_{j} \) in \( S \).
  
  Let \( \tau=\langle f_{1},..,f_{m}\rangle \).
  Consider the cone \( \delta =\langle e_{i} , f_{1},..,f_{m}\rangle \).
   If \( \delta \) is singular, 
  then, by Lemma \ref{miz}, \( v_{j} \) is not minimal in \( 
  S'=\bigcup_{\tau'<\delta:\operatorname{singular}}{\tau'}^o\cap N\), 
  therefore by \cite[Lemma 3.15]{i-k}, 
  we have a non-singular subdivision \( \Delta \) of \( \delta \) in which \( 
  \langle v_{j}\rangle \) 
  does not appear as a one-dimensional cone
   and every non-singular face of \( \delta \) does not 
  change.
  Here, \( \Delta \) is obtained by successive star-shaped 
  subdivision by centers \( \lambda_{1},\ldots,\lambda_{l} \).
  
  Now, we construct a subdivision of \( \sigma \). 
  
  {\bf Step 1.} Take the star-shaped subdivision \( \Sigma_{1} \) of  
  \( \sigma \) by 
  \( e_{i} \). 
  Here, the cone \( \delta \) appears in \( \Sigma_{1} \).
  Note that the morphism corresponding to this subdivision is 
  isomorphic outside of \( \overline{orb \tau_{i}} \).
  
  {\bf Step 2.} If \( \Sigma_{1} \) is simplicial, then put 
  \( \Sigma_{2}=\Sigma_{1} \).
  If \( \Sigma_{1} \) is not simplicial, then take a one-dimensional 
  face \( \langle \lambda \rangle \)  of  a minimal dimensional 
  non-simplicial cone and make the star-shaped subdivision of 
  \(\Sigma_{1} \) by \( 
  \lambda \).
  Then,  simplicial cones in \( \Sigma_{1} \) do  not change. 
  Continuing this procedure, we obtain a simplicial subdivision
  \( \Sigma_{2} \).
  
  {\bf Step 3.} If \( \delta \) is non-singular, then put \( 
  \Sigma_{3}=\Sigma_{2} \).
  If \( \delta \) is singular,
   take the successive star-shaped subdivisions \( \Sigma_{3} \) of 
  \( \Sigma_{2} \) with the centers \( \lambda_{1},\ldots,\lambda_{l} \).
   Then, the cone \( \delta \) in \( \Sigma_{2} \) is replaced by the fan 
  \( \Delta \).
  Note that this subdivision does not change non-singular cones of \( 
  \Sigma_{2} \), therefore the morphism corresponding to this 
  subdivision is isomorphic outside of the singularities.
  
  {\bf Step 4.} If \( \Sigma_{3} \) is singular,  take a cone \( \lambda=
  \langle p_{1},\ldots,p_{t}\rangle\in \Sigma_{2} \) with
  the maximal multiplicity.
  The  multiplicity is \( \vol P_{\lambda} \),
  where \( P_{\lambda}=\{\sum_{i=1}^t c_{i}p_{i}\mid 0\leq c_{i}<1\}\).
  Since \( \vol P_{\lambda}>1 \), there is a non-zero element 
  \(n'\in P_{\lambda}\cap 
  N \).
  Take the star-shaped subdivision with the center \( n' \).
  Then, the multiplicities of new cones on \( \lambda \) become less
  than \( \lambda \) and
  all  non-singular cones of \( \Sigma_{3} \) are unchanged.
  Continuing this procedure, we finally obtain a non-singular subdivision
  \( \Sigma_{4} \).  
  
  This subdivision \( \Sigma_{4} \) gives a birational morphism 
  \( \psi'':Y''\to \overline{X} \) which is isomorphic outside of 
  \( \overline{orb \tau_{i}}\cup \sing \overline{X} \).
  Therefore \( \nu\circ\psi'':Y''\to X \) is a resolution of 
  singularities of \( X \).
  As \( \langle v_{j}\rangle \) does not appear in \( \Sigma_{4} \) 
  as a cone, 
  the center of \( D_{v_{j}} \) on \( Y'' \) is contained in \( D_{e_{i}} \) 
  or some exceptional divisor \( D \) on \( Y'' \).
  Thus, \( D_{v_{j}} \) \( (j>w) \) is not an essential divisor over \( X \).
\end{pf}

\begin{thm}
  Let \( X \) be an affine pretoric variety. Then the Nash map
  \[ {\cal N}:\{{{Nash\ components\ in\ \pi_{X}^{-1}(\sing X)}}\} \to 
  \{{{essential\ divisors\ over\ X}}\} \]
  is bijective.
\end{thm}  

\begin{pf} 
  Let \( \varphi:Y \to X \) be a resolution in which 
   \( D_{e_{i}} \) \( (i=1,..,r) \) and \( D_{v_{j}} \) \( (j=1,.,w) \)
   appear.
  By Lemma \ref{mainlemma}, an essential divisor over \( X \) is 
  \( D_{e_{i}} \) \( (i=1,..,r) \) or \( D_{v_{j}} \) \( (j=1,.,w) \).
  Then, by the discussion in Example \ref{nash}, 
  the union of the Nash components is \[\bigcup_{i=1}^r 
  \overline{\varphi_{\infty}\pi_{Y}^{-1}(D_{e_{i}})}\cup
  \bigcup_{j=1}^w 
  \overline{\varphi_{\infty}\pi_{Y}^{-1}(D_{v_{j}})}.\]
  By Lemma \ref{equal} and Lemma \ref{bij}, it follows that
 \[  \overline{\varphi_{\infty}\pi_{Y}^{-1}(D_{e_{i}})}=
  \overline{T_{\infty}^X(e_{i})} \]
  \[ \overline{\varphi_{\infty}\pi_{Y}^{-1}(D_{v_{j}})}
  = \overline{T_{\infty}^X(v_{j})}. \]
  Here, each closed set does not contain any other set 
  by Lemma \ref{pretoricdominate},
  because each element of \(  \{{e_{i}}, 
   {v_{j}}\}^{i=1,.,r}_{{j=1,.,w}} \) is minimal among them.
  Hence, the number of the Nash components is \( r+w \), 
  while the number of the essential divisors is less than or equal to
  \( r+w \) by Lemma \ref{mainlemma}.
  Since the Nash map is injective, the both numbers must be \( r+w \) 
  and
   the Nash map is bijective.        
\end{pf}

\begin{cor}
  If \( X \) is a non-normal toric variety, then the Nash map for \( X 
  \) is bijective.
\end{cor}

For a general non-normal variety, we have a counter example to the Nash 
problem.
  It is obtained from the counter example in \cite{i-k}.
  
\begin{exmp}
  Let \( X \) be a non-normal variety with the normalization 
  \( \nu:\overline{X}\to X \) such that \( \overline{X} \) is the 
  hypersurface of \( \bC^5 \) defined by \( x_{1}^3+x_{2}^3+x_{3}^4
  +x_{4}^3+x_{5}^6=0 \).
  Then the Nash problem is negative for \( X \).
  
  Indeed, let \( \psi_{1}:Y_{1}\to \overline{X} \) be the blow-up at 
  \( 0 \), then \( Y_{1} \) has an isolated singularity \( P \).
  Then, let \( \psi_{2}:Y_{2}\to Y_{1} \) be the blow-up at \( P \).
  Let \( E_{1}\subset Y_{2} \) be the proper transform of the 
  exceptional divisor of \( \psi_{1} \) and \( E_{2} \) the exceptional
  divisor of \( \psi_{2} \).
  In \cite[\S 4]{i-k} we proved that \( E_{2} \) is not ruled and 
  \[ (\psi_{1}\psi_{2})_{\infty}(\pi_{Y_{2}}^{-1}(E_{2}))
  \subset
  \overline{(\psi_{1}\psi_{2})_{\infty}(\pi_{Y_{2}}^{-1}(E_{1}))}. \]
  Therefore, it follows that \( E_{2}  \) is an essential divisor 
  over \( X \) and 
  \[ (\nu\psi_{1}\psi_{2})_{\infty}(\pi_{Y_{2}}^{-1}(E_{2}))
  \subset
  \overline{(\nu\psi_{1}\psi_{2})_{\infty}(\pi_{Y_{2}}^{-1}(E_{1}))}, \]
  which shows that there  is no Nash component corresponding to 
  \( E_{2} \) (see \ref{nash})
  
  To construct such an example concretely, 
  we can define \( X \) as follows: 
  Let \( e_{1},..,e_{5} \) be a basis of \( M=\bZ^5 \) and \( \Gamma 
  \) the subsemigroup of \( M \) generated by \( e_{1},..,e_{4},
  2e_{5}, 3e_{5} \).
  Then, the canonical morphism 
  \( \nu:\bC^5\to  \spec \bC[\Gamma] \)
  induced from the injection 
  \( \Gamma \hookrightarrow \oplus_{i=1}^5 \bZ_{\geq 0}e_{i} \)
  is the normalization of a non-normal toric variety \( \spec \bC[\Gamma] \).
  It is sufficient to let \( X \) be the image \( \nu(\overline{X}) \).

\end{exmp}


\makeatletter \renewcommand{\@biblabel}[1]{\hfill#1.}\makeatother

\end{document}